\documentclass[a4paper, reqno]{amsart} 

\usepackage{amsmath}
\usepackage{amssymb}
\usepackage{color} 
\usepackage{graphicx} 
\usepackage{dsfont}
\usepackage[oldenum, olditem, alwaysadjust]{paralist}
 
\usepackage{hyperref} 
\hypersetup{ 
  pdftitle={A variational principle for weighted Delaunay triangulations and
    hyperideal polyhedra},
  pdfauthor={Boris A. Springborn}, 
  pdfsubject={},
  pdfkeywords={}, 
  pdfpagemode=None, 
  pdfstartview=FitH, 
  pdfnewwindow=true,
  breaklinks=true,
  pagebackref=true,
  colorlinks=true,
  linkcolor=black,
  anchorcolor=black,
  citecolor=black,
  filecolor=black,
  menucolor=black,
  pagecolor=black,
  urlcolor=black,
  bookmarksopen=false,
}

\graphicspath{{./figs/}}

\newcommand{\ie}{{\it i.e.}}

\newcommand{\etal}{{\it et al.}}

\newcommand{\R}{{\mathds R}}
\newcommand{\Rplus}{\R_{>0}}

\newcommand{\A}{{\mathcal A}}
\newcommand{\T}{{\mathcal T}}
\newcommand{\Vol}{\operatorname{Vol}}

\font\cyr=wncyr10 
\newcommand{\ML}{\operatorname{\mbox{\cyr L}}}

\newtheorem{theorem}{Theorem}
\newtheorem*{theorem*}{Theorem}
\newtheorem{lemma}{Lemma}
\newtheorem*{claim*}{Claim}
\theoremstyle{definition}
\newtheorem*{definition}{Definition}
\theoremstyle{remark}
\newtheorem*{remark}{Remark}

\title[Variational\! principle\! for\! weighted\! Delaunay\! triangulations]%
{A variational principle for weighted Delaunay triangulations and
  hyperideal polyhedra}

\author{Boris A.~Springborn}
\address{Boris Springborn\\
  Technische Universit\"at Berlin\\
  Institut f\"ur Mathematik, MA 8-3\\
  Strasse des 17.~Juni 136\\
  10623 Berlin, Germany\\}
\email{boris.springborn@tu-berlin.de} 
\thanks{Supported by the DFG Research Center \textsc{Matheon} in Berlin.}

\date{} 

\begin{document}

\begin{abstract}
  We use a variational principle to prove an existence and uniqueness theorem
  for planar weighted Delaunay triangulations (with non-intersect\-ing
  site-circles) with prescribed combinatorial type and circle intersection
  angles. Such weighted Delaunay triangulations may be interpreted as images
  of hyperbolic polyhedra with one vertex on and the remaining vertices
  beyond the infinite boundary of hyperbolic space. Thus the main theorem
  states necessary and sufficient conditions for the existence and uniqueness
  of such polyhedra with prescribed combinatorial type and dihedral angles.
  More generally, we consider weighted Delaunay triangulations in piecewise
  flat surfaces, allowing cone singularities with prescribed cone angles in
  the vertices. The material presented here extends work by Rivin on Delaunay
  triangulations and ideal polyhedra. 
\end{abstract}

\maketitle

\section{Introduction}

\subsection{Overview}

Rivin developed a variational method to prove the existence and uniqueness of
ideal hyperbolic polyhedra with prescribed combinatorial type and dihedral
angles, or equivalently, of planar Delaunay triangulations with prescribed
combinatorial type and circumcircle intersection angles~\cite{rivin94}. The
purpose of this article is to extend this method to hyperideal polyhedra and
to weighted Delaunay triangulations. Consider a finite set of disjoint
circular disks in the plane. For every triple of such disks there exists a
circle that intersects the boundaries of the disks orthogonally. The weighted
Delaunay triangulation induced by the disks consists of the triangles whose
vertices are the centers of a triple of disks such that the orthogonal circle
of this triple intersects no other disk more than orthogonally (see
Figure~\ref{fig:hyperdelaunay}). Such weighted Delaunay triangulations
correspond to hyperbolic polyhedra with one vertex on the infinite boundary
and all other vertices outside. The dihedral angles correspond to the
intersection angles of the orthogonal circles. More generally, we allow cone
singularities at the centers of the disks, and we call such weighted Delaunay
triangulation with cone singularities euclidean hyperideal circle patterns.
The question we consider is: Given an abstract triangulation $\T$ with vertex
set $V$, edge set $E$, and face set $T$, and given intersection angles
$\theta$ as a function on the edges and cone angles $\Xi$ as a function on
the vertices, does there exist a corresponding euclidean hyperideal circle
pattern, and is it unique? The main result is the following.

\begin{theorem}\label{thm:main}
  A euclidean hyperideal circle pattern with triangulation $\T$, intersection
  angles $\theta:E\rightarrow[0,\pi)$ and cone/boundary angles\/
  $\Xi:V\rightarrow(0,\infty)$ exists if and only if the set of coherent
  angle systems $\A(\T,\theta,\Xi)$ is not empty. In this case, the circle
  pattern is unique up to scale.
\end{theorem}

Section~\ref{sec:basics} contains the basic definitions and a precise
statement of the ``circle pattern problem'' under consideration. The set of
coherent angle systems $\A(\T,\theta,\Xi)$ is defined in
Section~\ref{sec:cas}. It is a subset of $\R^{6|T|}$ defined by a system of
linear equations and linear strict inequalities. The proof of
Theorem~\ref{thm:main} is based on the variational principle that is
presented in Section~\ref{sec:v_princ}, where a function
$F:\R^{6|T|}\rightarrow\R$ is defined explicitely in terms of Milnor's
Lobachevsky function.  The critical points of $F$ in $\A(\T,\theta,\Xi)$
correspond to solutions of the circle pattern problem (Lemma~\ref{lem:crit}).
The uniqueness of a solution follows immediately from the fact that the
function $F$ is strictly concave on $\A(\T,\theta,\Xi)$
(Lemma~\ref{lem:concave}).  To prove the existence of a solution, we show
that $F$ cannot attain its maximum on the boundary
$\overline{\A(\T,\theta,\Xi)}\setminus\A(\T,\theta,\Xi)$
(Lemma~\ref{lem:max}).
Sections~\ref{sec:proof_lem_crit}--\ref{sec:proof_lem_max} are devoted to the
proofs of these three main lemmas. The explicit formula for the function $F$
is based on a hyperbolic volume formula which is derived in
Section~\ref{sec:volumes}.

The variational principle presented here is not only a tool to prove the
existence and uniqueness Theorem~\ref{thm:main}. Since it reduces the circle
pattern problem to a convex optimization problem with linear constraints, it
provides a means for its numerical solution. This is important in view of
possible applications, such as using circle patterns to map 3D triangle
meshes to the plane.

Different necessary and sufficient conditions for the existence of hyperideal
circle patterns were obtained by Bao~\& Bonahon~\cite{bao02} and by
Schlenker~\cite{schlenker05b}. We discuss these results in
Section~\ref{sec:related_work}.

\subsection{Delaunay triangulations and hyperbolic polyhedra}
\label{sec:cp+hp}

\begin{figure}
  \centering
  \includegraphics[width=0.5\textwidth]{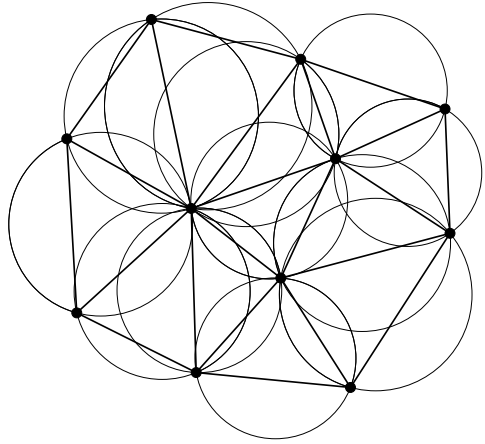}
  \caption{A Delaunay triangulation.}
  \label{fig:delaunay}
\end{figure}

``Patterns of circles'' have become objects of mathematical interest after
Thurston introduced them as elementary and intuitive images of polyhedra in
hyperbolic $3$-space~\cite{thurston}. Thus, a planar Delaunay
triangulation (\ie~a triangulation of a convex polygonal region with the
property that the circumcircle of each triangle does not contain any vertices
in its interior, see Figure~\ref{fig:delaunay})
can be viewed as representation of a convex hyperbolic polyhedron with all
vertices in the infinite boundary of hyperbolic space: Erase all interior
edges, keep only the circumcircles and the boundary edges, and extend the
boundary edges to straight lines. Consider the paper plane as the infinite
boundary of hyperbolic space, represented in the Poincar\'e half-space model.
In this model, hyperbolic planes are represented by hemispheres and
half-planes which intersect the infinite boundary orthogonally in circles and
lines, respectively. Therefore, if we erect hemispheres and orthogonal
half-planes over the circumcircles and the prolonged boundary edges, we
obtain a set of hyperbolic planes which bound a convex polyhedron. This
polyhedron's vertices are the vertices of the Delaunay triangulation and one
additional point, the infinite point of the boundary plane, where all the
hyperbolic planes corresponding to the boundary edges intersect. Moreover,
the dihedral angle at an edge of the polyhedron is equal to the angle in
which the corresponding circles/lines intersect.

A point in the infinite boundary of hyperbolic space is called an \emph{ideal
  point}, and a polyhedron with all vertices in the ideal boundary is called
an \emph{ideal polyhedron}. This terminology has become widely accepted; in
the old literature, the term ``ideal point'' was used for points
\emph{beyond} the ideal boundary, which are now called \emph{hyperideal
  points}. If we consider Delaunay triangulations up to similarity and
hyperbolic polyhedra up to isometry, then the construction above establishes
a $1$-to-$1$ correspondence between planar Delaunay triangulations and convex
ideal polyhedra with one marked vertex. (Essentially the same construction,
but represented in the projective model of hyperbolic space with a paraboloid
as the absolute quadric, is known in Discrete Geometry as the ``convex hull
construction''.)

\begin{figure}
  \centering
  \includegraphics[width=0.66\textwidth]{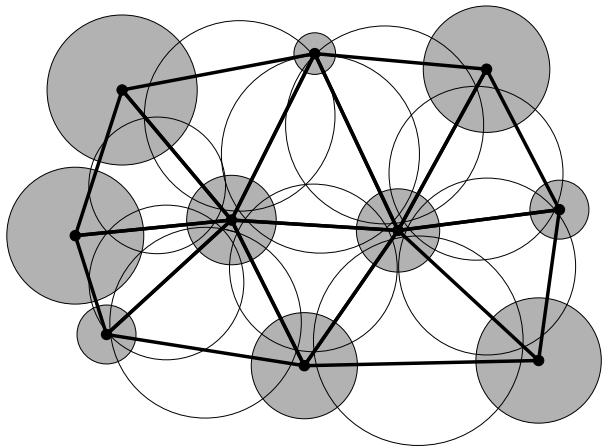}
  \caption{A weighted Delaunay triangulation. Note that the intersection
    points of two adjacent face-circles lie on the edge between them.}
  \label{fig:hyperdelaunay}
\end{figure}

\emph{Weighted Delaunay triangulations} are a well known generalization of
Delaunay triangulations~\cite{edelsbrunner01}, where the sites are not points
but circles (vertex-circles). We consider only the case where the closed
disks bounded by the vertex-circles are pairwise disjoint. Instead of an
empty circumcircle, to each triangle there corresponds a circle (face-circle)
which intersects the adjacent vertex-circles orthogonally and which
intersects no vertex-circle more than orthogonally (see
Figure~\ref{fig:hyperdelaunay}).
Weighted Delaunay triangulations with non-intersecting vertex circles
correspond to hyperbolic polyhedra with hyperideal vertices but with edges
still intersecting hyperbolic space. Hyperideal points are not represented as
such in the Poincar\'e half-space model. In the projective model, they are
simply represented by the points outside the absolute quadric. The plane that
is polar to such a point (with respect to the absolute quadric) intersects the
absolute quadric, hence it represents a hyperbolic plane. Thus, hyperideal
points are in $1$-to-$1$ correspondence with hyperbolic planes. Two
hyperbolic planes intersect orthogonally iff one is incident with the
hyperideal point polar to the other. The correspondence between Delaunay
triangulations and hyperbolic polyhedra extends to weighted Delaunay
triangulations: \emph{A weighted Delaunay triangulation with non-intersecting
  vertex circles corresponds to a convex hyperbolic polyhedron with one
  marked ideal vertex and all other vertices hyperideal, and with edges
  intersecting hyperbolic space.}  The problem we will consider is to
construct such polyhedra with prescribed combinatorial type and dihedral
angles. More generally, we will also consider weighted Delaunay
triangulations in piecewise flat surfaces with cone singularities with
prescribed cone angle at the vertices. These correspond to certain
non-compact hyperbolic cone manifolds with polyhedral boundary where the
lines of curvature connect one marked ideal vertex with the other hyperideal
vertices.

\subsection{Related work}
\label{sec:related_work}

For a comprehensive bibliography on circle packings and circle patterns we
refer to Stephenson's monograph~\cite{stephenson05}. Here, we can only attempt
to briefly discuss some of the the most important and most closely related
results. 

Let $\mathcal C$ be a cellulation of the 2-sphere and suppose there is a
weight $\theta_e\in(0,\pi)$ attached to each edge $e$.  Does there exist a
hyperbolic polyhedron that is combinatorially equivalent to $\mathcal C$ and
whose exterior dihedral angles are the weights $\theta_e$; and if so, is it
unique?  A complete answer to this question is unknown.  Andreev gave an
answer for compact polyhedra with non-obtuse dihedral
angles~\cite{andreev70a} (see also~\cite{roeder06}), and he extended his
result to non-obtuse angled polyhedra with finite volume, some or all
vertices of which may be in the sphere at infinity~\cite{andreev70b}. An
analogous existence and uniqueness theorem for circle patterns in surfaces of
non-positive Euler characteristic is due to Thurston~\cite{thurston}. The
intersection angles have to be non-obtuse but this theorem also allows the
circle pattern equivalent of hyperideal vertices. Chow~\& Luo~\cite{chow03}
gave a proof which is inspired by work on the Ricci flow on surfaces. They
also show that there is a variational principle for this circle pattern
problem, but only the derivatives of the functional are known.

Rivin classified convex ideal polyhedra without any restriction to
non-obtuse dihedral angles:

\begin{theorem}[Rivin~\cite{rivin96}]
  \label{thm:rivin96}
  There exists an ideal polyhedron that is combinatorially equivalent to a
  cellulation $\mathcal C$ of $S^2$ and that has prescribed exterior dihedral
  angles $\theta_e\in(0,\pi)$ iff for each cycle $\gamma$ in the 1-skeleton
  of the dual cellulation $\mathcal C^*$ the inequality
  \begin{equation}
    \label{eq:rivin_condition}
    \sum_{e^* \in\gamma} \theta_e \geq 2\pi
  \end{equation}
  holds and equality holds iff $\gamma$ is the boundary of a face
  of $\mathcal C^*$.  If it exists, the polyhedron is unique.
\end{theorem}

Bao and Bonahon generalized Rivin's result for polyhedra with all vertices
\emph{on or beyond} the sphere at infinity:

\begin{theorem}[Bao \& Bonahon~\cite{bao02}]
  \label{thm:bao}
  There exists a polyhedron with ideal and hyperideal vertices and with edges
  intersecting hyperbolic space that is combinatorially equivalent to a
  cellulation $\mathcal C$ of $S^2$ and has prescribed exterior dihedral
  angles $\theta_e\in(0,\pi)$ iff the following conditions hold:
  
  (i) For each cycle $\gamma$ in the 1-skeleton of the dual cellulation
  $\mathcal C^*$ the inequality~\eqref{eq:rivin_condition}
  holds and equality may hold only if $\gamma$ is the boundary of a face
  of $\mathcal C^*$.

  (ii) For each simple path $\gamma$ in the 1-skeleton of the dual
  cellulation $\mathcal C^*$ that joins two different vertices of a face
  $v^*$ of $\mathcal C^*$ and that is not contained in the boundary of $v^*$,
  $\sum_{e^*\in\gamma}\theta_e>\pi$.

  If it exists, the polyhedron is unique.
\end{theorem}

Andreev, Rivin, and Bao~\& Bonahon obtain their results by employing variants
the \emph{method of continuity} (or \emph{deformation method}) that was
pioneered by Alexandrov~\cite{alexandrov05}. Schlenker gave a different proof
of Bao~\& Bonahon's Theorem~\cite{schlenker05a}.

A variational approach to construct ideal hyperbolic polyhedra and, more
generally, Delaunay triangulations of piecewise flat surfaces was also
provided by Rivin~\cite{rivin94}. The basic idea is to build an ideal
polyhedron by gluing together ideal tetrahedra, or equivalently, to build a
Delaunay triangulation by gluing together triangles.  The angles of the
triangles are considered as variables. They have to satisfy simple linear
equality and inequality constraints: They have to be positive and the three
angles in each triangle have to sum to $\pi$. The angles around a vertex have
to sum to $2\pi$ (more generally, some specified cone angle).  Finally, to
get the right circle intersection angles, the angles opposite an edge $e$
have to sum to $\pi-\theta_e$. Using Colin de~Verdi\`ere's
terminology~\cite{colin91}, we call an assignment of values to the angle
variables that satisfies these constraints a \emph{coherent angle system}. A
coherent angle system does in general not represent a Delaunay triangulation
because for the triangles to fit together, further non-linear conditions on
the angles have to be satisfied. Nevertheless, the following theorem holds.

\begin{theorem}[Rivin~\cite{rivin94}]
  \label{thm:rivin94}
  Let $\mathcal S$ be a cellulation of a surface, let an intersection angle
  $\theta_e\in[0,\pi)$ be assigned to each edge $e$ and a cone angle $\Xi_v$
  be assigned to each vertex $v$. There exists a Delaunay triangulation of a
  piecewise flat surface that is combinatorially equivalent to $\mathcal S$
  and has intersection angles $\theta$ and cone angles $\Xi$ if and only if
  the above constraints on the angle variables are feasible, \ie~if a
  coherent angle system exists. In that case, the Delaunay triangulation is
  unique up to scale.
\end{theorem}

Note that the necessary and sufficient conditions of
Theorem~\ref{thm:rivin96}, involving inequalities of sums of intersection
angles over paths, are very different from the conditions of
Theorem~\ref{thm:rivin94}, involving the existence of a coherent angle
system. There is also a third type of conditions for the existence of a
Delaunay triangulation of a piecewise flat surface with prescribed
intersection angles and cone angles that was first obtained by
Bowditch~\cite{bowditch91}. It is by no means a triviality to directly derive
one type of conditions from another type~\cite{rivin03} \cite{bobenko04}. A
variant of Rivin's variational approach for Delaunay decompositions of
\emph{hyperbolic} surfaces was developed by Leibon~\cite{leibon02}. In this
article, we extend Rivin's variational approach to euclidean weighted
Delaunay triangulations with non-intersecting vertex circles. The main
Theorem~\ref{thm:main} is of the type ``a weighted Delaunay triangulation
exists uniquely iff a coherent angle system exists.'' The scope of
Theorem~\ref{thm:main} has non-empty intersection with Bao~\& Bonahon's
Theorem~\ref{thm:bao}. Both cover hyperbolic polyhedra with hyperideal
vertices and precisely one ideal vertex.  But the conditions are of a
different type and not obviously equivalent.  Theorem~\ref{thm:main} also
covers weighted Delaunay triangulations in flat tori and, more generally, in
piecewise flat surfaces, possibly with boundary.  The variational method is
better suited for numerical computation.

Recently, Schlenker has treated weighted Delaunay triangulations in piecewise
flat and piecewise hyperbolic surfaces using a deformation
method~\cite{schlenker05b}. He obtains an existence and uniqueness theorem
\cite[Theorem~1.4]{schlenker05b} with the same scope as
Theorem~\ref{thm:main}, but the conditions for existence are in terms of
angle sums over paths like in Theorem~\ref{thm:bao}. This seems to be the
first time that this type of conditions was obtained for circle patterns with
cone singularities. It would be interesting to show directly that the
conditions of his theorem are equivalent to the conditions of
Theorem~\ref{thm:main}.

Circle patterns have been applied to map 3D triangle meshes quasi-con\-for\-mal\-ly
to the plane. To improve these methods was an important motivation for this
work. The group around Stephenson first used circle packings (circles
touching without overlap) to construct planar maps of the human
cerebellum~\cite{hurdal99}. This original method only takes the combinatorics
of the input mesh into account. A later version uses so called
\emph{inversive distance packings}~\cite{bowers03}. (The \emph{inversive
  distance} of two non-intersecting circles is the $\cosh$ of the hyperbolic
distance of the two planes they represent.) Inversive distance packings are
similar to the weighted Delaunay triangulations considered here except that
the inversive distances of the vertex-circles are prescribed instead of the
intersection angles of the face-circles, and there is no Delaunay criterion.
Unfortunately, no existence and uniqueness theorem for inversive distance
packings is known. Kharevych~\etal~\cite{kharevych06} proceed along a
different path.  They first read off the intersection angles between
circumcircles of the 3D triangle mesh. Then they construct a planar Delaunay
triangulation with intersection angles as close to the measured angles as
possible. To construct the Delaunay triangulation they use a variational
principle by Bobenko~\& Springborn~\cite{bobenko04}, which is related to
Rivin's via a Legendre transformation. It has the advantage that the
variables are (logarithmic) circle radii which are not subject to any
constraints. However, the prior angle-adjustment is achieved by solving a
quadratic programming problem, which is the most complicated and
computationally most expensive stage of the algorithm. One may hope that
using weighted Delaunay triangulations will provide a way to escape the tight
constraints that have to be satisfied by the intersection angles of a
Delaunay triangulation without giving up all mathematical certainty regarding
existence and uniqueness. Another interesting question is this: Can one
formulate a dual variational principle for weighted Delaunay triangulations,
with an \emph{explicit formula} for the functional, where the variables are
circle radii and inversive distances?

\section*{Acknowledgments}

\noindent
The research for this article was conducted almost entirely while I enjoyed
the hospitality of the \textit{Mathematisches Forschungsinstitut
  Oberwolfach}, where I participated in the Research in Pairs Program
together with Jean-Marc Schlenker, who was working on his closely related
paper~\cite{schlenker05b}. I am grateful for the excellent working conditions
I experienced in Oberwolfach and for the extremely inspiring and fruitful
discussions with Jean-Marc, who was closely involved in the work presented
here.

\section{Euclidean hyperideal circle patterns}
\label{sec:basics}

\subsection{Basic definitions}

A \emph{surface} is a two-di\-men\-sio\-nal manifold, possibly with boundary. A
\emph{triangulated surface} $\T$ (or \emph{triangulation} for short) is a
two-di\-men\-sio\-nal CW complex whose total space is a surface $S$ and which has
the property that for each two-cell attaching map $\sigma:B^2\rightarrow S$
the set $\sigma^{-1}(V)$ contains three points, where $V$ is the vertex set
(zero-skeleton) of the CW complex.

This definition allows non-regular triangulations. A cell complex is called
\emph{regular} if the cell attaching homomorphisms embed the closed cells. A
cell complex is called \emph{strongly regular}\/ if it is regular and the
intersection of two closed cells is empty or a closed cell. The usual
definition of simplicial complexes implies that they are strongly regular.
Throughout this paper we assume all triangulations to be regular, \emph{but
  only to simplify notation}. We will label vertices by $i$, $j$, $k$,
$\ldots$ and denote edges by pairs $ij$ and triangles by triples $ijk$.
However, this regularity assumption is not essential for the material
presented here; everything holds for non-regular triangulations as well.

We denote the set of vertices, edges, and triangles of a triangulation
$\T$ by $V$, $E$, and $T$, respectively.  Throughout this paper, all
triangulations are assumed to be finite.

A \emph{triangulated piecewise flat surface} $(\T, d)$ is a
triangulated surface $\T$ equip\-ped with a metric $d$ such that for
each two-cell attaching map $\sigma:B^2\rightarrow S$ the closed disk $B^2$
equipped with the pulled back metric $\sigma^* d$ is a euclidean triangle,
and $\sigma$ maps the vertices of this triangle to vertices of the CW
complex. In other words, a triangulated piecewise flat surface is a surface
obtained by glueing together euclidean triangles along their sides. (Of
course, sides that are identified by glueing must have the same length.) The
metric $d$ is flat except in the vertices of the triangulation where it may
have cone-like singularities. The \emph{cone angle}\/ at a vertex $i$ is the
sum of all triangle angles incident at $i$. If the cone angle at a vertex is
$2\pi$ then the metric is flat there. A triangulated piecewise flat surface
is determined by the triangulation $\T$ and the function
$l:E\rightarrow\Rplus$ that maps each edge $ij\in E$ to its length $l_{ij}$.
For each triangle $ijk\in T$, the lengths $l_{ij}$, $l_{jk}$, $l_{ki}$
satisfy the triangle inequalities. Conversely, a triangulation $\T$
and a function $l:E\rightarrow\Rplus$ that satisfies the triangle
inequalities for each triangle determines a triangulated piecewise flat
surface.

A \emph{euclidean hyperideal circle pattern} is a triangulated piecewise flat
surface together with a function $r:V\rightarrow \Rplus$ with the following
two properties.

\begin{compactitem}[(i)]
\item[(i)] For each edge $ij\in E$,
  $
  r_i + r_j < l_{ij},
  $
  where $l_{ij}$ is the length of the edge.
\end{compactitem}

\noindent%
Let $ijk\in T$ be a triangle of the triangulation $\T$. If we draw a
triangle with sides $l_{ij}$, $l_{jk}$ and $l_{ki}$ in the euclidean plane
and circles with radii $r_i$, $r_j$ and $r_k$ around the vertices, then the
property (i) simply says that these circles do not touch or intersect.
Consequently there exists a unique fourth circle that intersects all three
circles orthogonally. The second condition concerns these orthogonally
intersecting circles.

\begin{compactitem}[(ii)]
\item[(ii)] Let $ij\in E$ be an \emph{interior}\/ edge. Let $ijk$ and $jil$
  be the adjacent triangles on either side. (These may actually be one and
  the same triangle if the triangulation is not regular.) Draw two abutting
  triangles with the same side lengths in the euclidean plane, and draw
  circles with radii $r_i$, $r_j$, $r_k$ and $r_l$ around the vertices. Then
  the orthogonal circle through the vertex-circles of one triangle intersects
  the fourth vertex-circle either not at all or at an angle that is less
  than~$\frac{\pi}{2}$.
\end{compactitem}

In other words, a euclidean hyperideal circle pattern is a weighted Delaunay
triangulation with non-intersecting vertex-circles in a piecewise flat
surface. (Note that condition (ii) invokes an edge-local Delaunay condition.
This raises the question whether it is also true for piecewise flat surfaces
that the local condition implies the global condition that no face-circle
intersects any vertex-circle more than orthogonally. Such questions shall not
be treated here. We refer to Bobenko~\& Springborn~\cite{bobenko05} and
Bobenko~\& Izmestiev~\cite{bobenko} for a more thorough treatment of
Delaunay triangulations and weighted Delaunay triangulations in piecewise
flat surfaces.)

\begin{figure}
  \centering \input{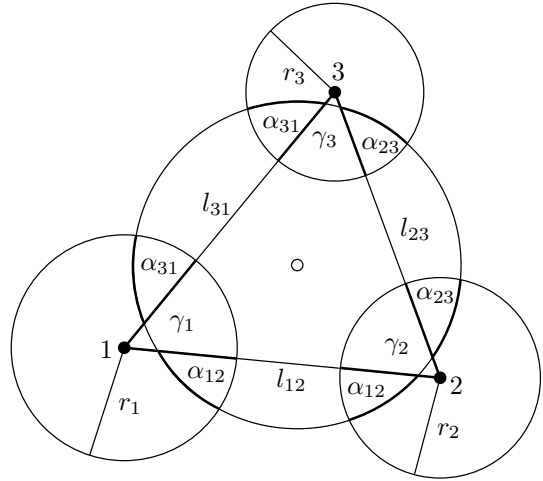}
  \caption{Local geometry at a triangle.} 
  \label{fig:triangle}
\end{figure}

Just as a triangulated piecewise flat surface is obtained by glueing together
euclidean triangles along the edges, a euclidean hyperbolic circle pattern is
obtained by putting together triangles which are decorated by circles as
shown in Figure~\ref{fig:triangle}.

\subsection{Interpretation as hyperbolic polyhedra}
\label{sec:hyperbolic_polyhedra}

\begin{figure}
  \centering
  \input{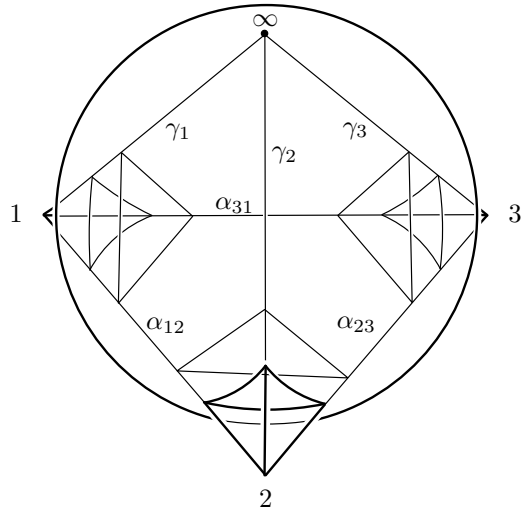}
  \caption{A tetrahedron in hyperbolic space (shown in the projective model)
    with one ideal vertex and three hyperideal vertices. The tetrahedron is
    truncated by the polar planes of the hyperideal vertices.} 
  \label{fig:truncated_tet}
\end{figure}

By the construction explained in Section~\ref{sec:cp+hp}, the decorated
triangle shown in Figure~\ref{fig:triangle} may be interpreted as a
tetrahedron in three-di\-men\-sio\-nal hyperbolic space with one vertex on the
sphere at infinity and three vertices beyond that sphere, as shown in
Figure~\ref{fig:truncated_tet}.
The sides of the triangle and the face-circle that intersects the three
vertex-circles orthogonally correspond to hyperbolic planes that bound a
tetrahedron with one ideal and three hyperideal vertices, which are
represented by the vertex-circles.

If we can put together the decorated triangles to form a circle pattern, we
can also glue the corresponding hyperbolic tetrahedra to form a hyperbolic
cone manifold with polyhedral boundary. Such a glueing will identify the ideal
vertices of all tetrahedra. The resulting point will either be an ideal
vertex if $\T$ has boundary or a cusp of the hyperbolic manifold if $\T$ is
closed. The hyperideal vertices of the tetrahedra will be identified in
groups corresponding to vertices of $\T$ to form hyperideal vertices of the
polyhedral boundary. There will in general be cone lines running from the cusp
or ideal vertex to the hyperideal vertices. If the circle pattern has no
curvature at the vertices, we obtain a non-compact hyperbolic manifold with
polyhedral boundary. If in addition the triangulation $\T$ is topologically a
disk, we obtain a convex hyperbolic polyhedron with one ideal vertex all
others hyperideal.

\subsection{The circle pattern problem}
\label{sec:problem}

From a euclidean hyperideal circle pattern one can read off the following
data:

\smallskip
\begin{compactitem}
\item A triangulated surface $\T$.
\item For each vertex $i\in V$, the sum $\Xi_i\in(0,\infty)$ of incident
  triangle angles. For an interior vertex, this is the cone angle at
  $i$, and $\Xi_i=2\pi$ if the circle pattern is flat at $i$. For a boundary
  vertex, $\Xi_i$ is the interior angle of the polygonal boundary at $i$.
\item For each interior edge $ij\in E$, the intersection angle
  $\theta_{ij}\in[0,\pi)$ of the orthogonal circles corresponding to the two
  adjacent triangles $ijk$ and $jil$. An intersection angle $\theta_{ij}=0$
  means that the orthogonal circles of triangles $ijk$ and $jil$ coincide.
\item For each boundary edge $ij\in E$, the intersection angle, also
  denoted by $\theta_{ij}$, of the orthogonal circle corresponding to the
  adjacent triangle $ijk$ with the line segment containing the edge
  $ij$. The range of these intersection angles at boundary edges is
  $\theta_{ij}\in(0,\pi)$. 
\end{compactitem}

\smallskip
We consider the following \emph{circle pattern problem}: \emph{Given}\/ an
abstract triangulation $\T$ and angle data $\Xi:V\rightarrow(0,\infty)$,
$\theta:E\rightarrow[0,\pi)$ \emph{find} a corresponding euclidean hyperideal
circle pattern.

\subsection{Local geometry at a triangle}
\label{sec:local_geom}

Consider a geometric figure consisting of a euclidean triangle with three
non-touching and non-inter\-sec\-ting circles centered at the vertices and a
fourth circle intersecting the other three orthogonally. Let
$\alpha_{12}$, $\alpha_{23}$, $\alpha_{31}$, $\gamma_1$, $\gamma_2$, $\gamma_3$
be the angles shown in Figure~\ref{fig:triangle}.
They are positive,
\begin{equation}
  \label{eq:alpha_gamma_positive}
      \alpha_{ij}>0,\quad\gamma_j>0,
\end{equation}
satisfy the angle sum equation
\begin{equation}
  \label{eq:gamma_triangle_sum}
      \gamma_1+\gamma_2+\gamma_3=\pi,  
\end{equation}
and the inequalities
\begin{equation}
  \label{eq:alpha_gamma_inequalities}
  \begin{split}
      \gamma_1+\alpha_{12}+\alpha_{31}&<\pi, \\
      \gamma_2+\alpha_{23}+\alpha_{12}&<\pi, \\
      \gamma_3+\alpha_{31}+\alpha_{23}&<\pi.    
  \end{split}
\end{equation}
Let
\begin{equation*}
  \Delta=
  \left\{(\alpha_{12},\alpha_{23},\alpha_{31},\gamma_1,\gamma_2,\gamma_3)
    \in\R^6
    \text{ satisfying \eqref{eq:alpha_gamma_positive}, 
      \eqref{eq:gamma_triangle_sum}, and 
      \eqref{eq:alpha_gamma_inequalities}}\right\}.
\end{equation*}

Conversely, if
$(\alpha_{12},\alpha_{23},\alpha_{31},\gamma_1,\gamma_2,\gamma_3)\in\Delta$,
then there exists one such figure with these angles, and only one up
to similarity. Indeed, the construction of such a figure is simple:
Draw any circle in the plane. (This fixes the scale and translational
degrees of freedom.) Then draw a line intersecting it at the angle
$\alpha_{12}$. (This fixes the remaining rotational degree of
freedom.) Then draw the other two lines intersecting the circle and
the first line at the prescribed angles. The
inequalities~\eqref{eq:alpha_gamma_inequalities} ensure that the lines
intersect outside the face-circle and that the orthogonal
vertex-circles do not intersect each other.

\subsection{Coherent angle systems}
\label{sec:cas}

Let $\T$ be a triangulation with triangle set $T$, edge set $E$, and
vertex set $V$. We label the coordinates of points in $\R^{6|T|}$ by
$\alpha^t_{ij}$, $\gamma^t_i$ where $t\in T$ is a triangle and $i, j\in t$
are vertices of $t$. We fix this labeling once and for all. Let
$\Xi:V\rightarrow(0,\infty)$ be a function on the vertices and
$\theta:E\rightarrow[0,\pi)$ be a function on the edges. The \emph{space of
  coherent angle systems} $\A(\T,\theta,\Xi)$ is the set of all points
$(\alpha,\gamma)\in\R^{6|T|}$ such that

\smallskip
\begin{compactitem}
\item for each triangle $t=ijk\in T$
  \begin{equation*}
    (\alpha_{ij}^t,\alpha_{jk}^t,\alpha_{ki}^t,
    \gamma_i^t,\gamma_j^t,\gamma_k^t)
    \in \Delta,
  \end{equation*}
\item for each interior edge $ij\in E$
  \begin{equation}
    \label{eq:alpha_sum}
    \alpha_{ij}^{t}+\alpha_{ji}^{t'}=\pi-\theta_{ij},
  \end{equation}
  where $t,t'\in T$ are the adjacent triangles on either side of edge $ij$,
\item for each boundary edge $ij\in E$, $\theta_{ij}>0$ and
  \begin{equation}
    \label{eq:alpha_bdy}
    \alpha_{ij}^{t}=\pi-\theta_{ij},
  \end{equation}
  where $t$ is the triangle incident with edge $ij$,
\item for each vertex $i\in V$
  \begin{equation}
    \label{eq:gamma_vertex_sum}
    \sum_{t\in T:i\in t} \gamma^t_i=\Xi_i. 
  \end{equation}
\end{compactitem}

\smallskip
A \emph{coherent angle system}\/ is an element of $\A(\T,\theta,\Xi)$. If
$\A(\T,\theta,\Xi)$ is not empty, then the closure
$\overline{\A(\T,\theta,\Xi)}$ is a compact polytope in $\R^{6|T|}$, and
$\A(\T,\theta,\Xi)$ is its relative interior.

\subsection{Existence and uniqueness}
\label{sec:e+u}

The main Theorem~\ref{thm:main} reduces the question of existence and
uniqueness of a solution of the circle pattern problem of
Section~\ref{sec:problem} to a linear feasibility problem.
The ``only if'' part of the theorem---if a circle pattern exists then a
coherent angle system exists---is trivial, since one can simply read off a
coherent angle system from a euclidean hyperideal circle pattern. The ``if''
part---if a coherent angle system exists then a circle pattern exists---is
not trivial, and the rest of this paper is devoted to proving it.  It is
\emph{not} true that for each coherent angle system there exists a euclidean
hyperideal circle pattern with these angles. While it \emph{is} true that a
coherent angle system determines up to similarity a geometric figure as shown
in Figure~\ref{fig:triangle} for each triangle of the triangulation, these
figures cannot in general be put together to form a circle pattern. The
relative scale of the two figures corresponding to neighboring triangles is
determined by the condition that the triangle edges to be glued together must
have the same length. But the two pairs of corresponding vertex circles in
the two figures will in general not have matching radii. Even if we
disregarded the vertex circles, it is in general not be possible to choose
consistently a scale for each triangle such that the corresponding sides of
triangles that are to be glued together have the same length.

\section{A variational principle}
\label{sec:v_princ}

\noindent
The space $\A(\T,\theta,\Xi)\subset\R^{6|T|}$ of coherent angle systems is
defined by \emph{linear} equations and inequalities. For a coherent angle
system to describe a solution for the circle pattern problem, it has to
satisfy in addition certain \emph{non-linear} equations, which guarantee that
the triangle figures can be put together to form a circle pattern. It turns
out that these non-linear compatibility conditions are the conditions for a
critical point of the function $F:\R^{6|T|}\rightarrow\R$ defined below under
variations in $\A(\T,\theta,\Xi)\subset\R^{6|T|}$. This is the content of
Lemma~\ref{lem:crit}. Together with Lemmas~\ref{lem:concave}
and~\ref{lem:max} this provides a proof of Theorem~\ref{thm:main}.

\subsection{The functional} 
\begin{figure}
  \centering
  \includegraphics[width=0.5\textwidth]{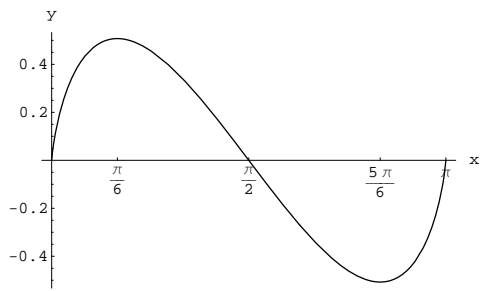}
  \caption{Milnor's Lobachevsky function, $y=\ML(x)$.}
  \label{fig:ML_plot}
\end{figure}

Label $\R^{6|T|}$ as in Section~\ref{sec:cas} and define
\begin{equation*}
  F_{\T}:\R^{6|T|}\rightarrow\R
\end{equation*}
by
\begin{equation}
  \label{eq:F_sum}
  F_{\T}=\sum_{t=ijk\in T} V(\alpha_{ij}^t,\alpha_{jk}^t,\alpha_{ki}^t,
    \gamma_i^t,\gamma_j^t,\gamma_k^t),
\end{equation}
where
\begin{equation}
  \label{eq:V}
  \begin{split}
    &2V(\alpha_{12},\alpha_{23},\alpha_{31},\gamma_1,\gamma_2,\gamma_3) = \\
  &\begin{alignedat}{3}
    &\hspace{4em}\ML(\gamma_1) 
    &+&\hspace{4em}\ML(\gamma_2) 
    &+&\hspace{4em}\ML(\gamma_3)    \\
    +&\ML\big(\tfrac{\pi+\alpha_{31}-\alpha_{12}-\gamma_1}{2}\big)
    &+&\ML\big(\tfrac{\pi+\alpha_{12}-\alpha_{23}-\gamma_2}{2}\big)
    &+&\ML\big(\tfrac{\pi+\alpha_{23}-\alpha_{31}-\gamma_3}{2}\big)\\
    +&\ML\big(\tfrac{\pi-\alpha_{31}+\alpha_{12}-\gamma_1}{2}\big)
    &+&\ML\big(\tfrac{\pi-\alpha_{12}+\alpha_{23}-\gamma_2}{2}\big)
    &+&\ML\big(\tfrac{\pi-\alpha_{23}+\alpha_{31}-\gamma_3}{2}\big)\\
    +&\ML\big(\tfrac{\pi+\alpha_{31}+\alpha_{12}-\gamma_1}{2}\big)
    &+&\ML\big(\tfrac{\pi+\alpha_{12}+\alpha_{23}-\gamma_2}{2}\big)
    &+&\ML\big(\tfrac{\pi+\alpha_{23}+\alpha_{31}-\gamma_3}{2}\big)\\
    +&\ML\big(\tfrac{\pi-\alpha_{31}-\alpha_{12}-\gamma_1}{2}\big)
    &+&\ML\big(\tfrac{\pi-\alpha_{12}-\alpha_{23}-\gamma_2}{2}\big)
    &+&\ML\big(\tfrac{\pi-\alpha_{23}-\alpha_{31}-\gamma_3}{2}\big),
  \end{alignedat}
\end{split}
\end{equation}
and the function 
\begin{equation*}
  \ML(x)=-\int_0^x\log|2\sin\xi|\,d\xi
\end{equation*}
is \emph{Milnor's Lobachevsky function}\/~\cite{milnor82},
\cite{milnor94:_volume_in_hyp}. (This is up to scale the same as
\emph{Clausen's integral} $\operatorname{Cl_2}(x)=2\ML(\frac{x}{2})$; see
Clausen~\cite{clausen32}, Lewin~\cite{lewin81}.) The function $\ML$ is
$\pi$-periodic, continuous, and odd. It is smooth everywhere except at
integer multiples of $\pi$ where its graph has a vertical tangent; see
Figure~\ref{fig:ML_plot}.
We will simply write $F$ for $F_{\T}$ when the triangulation $\T$ can be
inferred from the context.

\subsection{The main lemmas}
\label{sec:main_lemmas}

The following three lemmas imply Theorem~\ref{thm:main}. By
Lemma~\ref{lem:crit}, the critical points of $F$ in $\A(\T,\theta,\Xi)$
correspond to the solutions of the circle pattern
problem. Lemma~\ref{lem:concave} implies the uniqueness claim of
Theorem~\ref{thm:main}. The existence claim follows from Lemma~\ref{lem:max}.

\begin{lemma}
  \label{lem:crit}
  A coherent angle system $p\in\A(\T,\theta,\Xi)$ is a critical point of~$F$
  under variations in $\A(\T,\theta,\Xi)$ if and only if the decorated
  triangles fit together.
\end{lemma}

\begin{lemma}
  \label{lem:concave}
  The function $F$ is strictly concave on $\A(\T,\theta,\Xi)$. 
\end{lemma}

\begin{lemma}
  \label{lem:max}
  If $\A(\T,\theta,\Xi)$ is non-empty, then the restriction of $F$ to the
  closure $\overline{\A(\T,\theta,\Xi)}$ attains its maximum in
  $\A(\T,\theta,\Xi)$.
\end{lemma}

\section{Proof of Lemma~\ref{lem:crit}}
\label{sec:proof_lem_crit}

\noindent
The key ingredient to this proof is Schl\"afli's formula for the derivative
of the volume of a hyperbolic polyhedron when it is deformed in such a way
that its combinatorial type is preserved; see Milnor~\cite{milnor94},
\cite{schlaefli50}.

\begin{theorem*}[Schl\"afli's differential volume formula]
  The differential of the volume function $V$ on the space of 3-di\-men\-sio\-nal
  hyperbolic polyhedra of a fixed combinatorial type is
  \begin{equation}
    \label{eq:schlaefli}
    dV=-\frac{1}{2}\sum_{ij} a_{ij}\,d\alpha_{ij},
  \end{equation}
  where the sum is taken over the edges $ij$, and $a_{ij}$, $\alpha_i$ are
  the length and interior dihedral angle at edge $ij$.
\end{theorem*}

Hyperbolic polyhedra with some or all vertices on the sphere at infinity
still have finite volume. Milnor notes (in the concluding remarks
of~\cite{milnor94}) that Equation~\eqref{eq:schlaefli} remains true for such
polyhedra, under the following modification which is necessary because the
edges incident with an ideal vertex are of course infinitely long. Choose
arbitrary horospheres centered at the ideal vertices, and for edges $ij$
incident with an ideal vertex let $a_{ij}$ be the length of the edge
truncated at the horosphere(s) centered the ideal endpoint(s). One should add
that only deformations that leave the ideal vertices on the infinite sphere
are considered. The angle sum at an ideal vertex remains constant under such
a deformation, so $\sum d\alpha_{ij}=0$. If one chooses a different
horosphere at an ideal vertex $i$, the resulting truncated edge lengths of
the incident edges differ by the same additive constant: $\alpha_{ij}$
becomes $\alpha_{ij}+c_i$. Hence the right hand side of
Equation~\eqref{eq:schlaefli} does not depend on the choice of
horospheres. A proof for this extension of Schl\"afli's differential volume
formula to polyhedra with ideal vertices is contained in~\cite{springborn03}
(Lemma~4.1).

Next we extend Schl\"afli's differential volume formula to polyhedra that
have vertices beyond the sphere at infinity, but all of whose edges still
intersect hyperbolic space. If we truncate the hyperideal vertices at the
dual hyperbolic planes, we obtain a polyhedron with finite volume. (It may
still have ideal vertices.)

\begin{definition}
  \label{pag:def_trunc_vol}
  The \emph{truncated volume} of a hyperbolic polyhedron with vertices beyond
  infinity is the finite volume of the corresponding truncated polyhedron,
  truncated at the hyperbolic planes dual to the hyperideal vertices.
\end{definition}

Equation~\eqref{eq:schlaefli} remains true for polyhedra with hyperideal
vertices if we let $V$ be the truncated volume and $a_{ij}$ be the edge
lengths truncated at the dual planes of hyperideal vertices and at
horospheres centered at ideal vertices. Indeed, if we apply Schl\"aflis
differential volume formula to the truncated polyhedron, the edges introduced
by the truncation (those between original faces and truncation planes) have
dihedral angle $\pi/2$, and this angle is constant during any deformation. So
for these edges $d\alpha_{ij}=0$. Thus, only terms involving the original
edges remain in~\eqref{eq:schlaefli}.

\begin{lemma}
  \label{lem:volume_V}
  The truncated volume of the tetrahedron with one ideal and three
  hyperideal vertices with angles as shown in
  Figure~\ref{fig:truncated_tet} is equal to
  $V(\alpha_{12},\alpha_{23},\alpha_{31},\gamma_1,\gamma_2,\gamma_3)$
  as defined by Equation~\eqref{eq:V}.
\end{lemma}

\noindent
We prove this lemma in Section~\ref{sec:volumes}. 

\begin{remark}
  \label{pag:rem_same_formula}
  Exactly the same formula holds for tetrahedra with one ideal and three
  \emph{finite} vertices~\cite{vinberg91}, \cite{vinberg93}. The only
  difference is that the dihedral angles for such polyhedra satisfy the
  opposite inequalities instead of~\eqref{eq:alpha_gamma_inequalities}. This
  seems to be a common phenomenon: that the same volume formula holds
  regardless of whether certain vertices are finite or hyperideal. In the
  case of orthoschemes (the generalization of triply orthogonal tetrahedra to
  arbitrary dimension) this was observed by Kellerhals~\cite{kellerhals89}.
  Section~\ref{sec:volumes} contains some more examples. It would obviously
  be useful to have this observation cast into a theorem. The author is not
  aware that this has been done.
\end{remark}

From Schl\"afli's differential volume formula and Lemma~\ref{lem:volume_V} we
obtain the following.

\begin{lemma}
  \label{lem:V_derivative}
  If we choose a horosphere centered at the ideal vertex of the truncated
  hyperbolic tetrahedron shown in Figure~\ref{fig:truncated_tet} then 
  \begin{equation}
    \label{eq:a_ij}
    -2 \frac{\partial V}{\partial\alpha_{ij}} = a_{ij},
  \end{equation}
  where $a_{ij}$ is the length of the edge between the hyperideal vertices
  $i$ and $j$ truncated at the polar planes, and
  \begin{equation}
    \label{eq:a_i}
    -2 \Big(\frac{\partial V}{\partial\gamma_{i}} - \frac{\partial
    V}{\partial\gamma_{j}}\Big) = a_{i} - a_{j},
  \end{equation}
  where $a_{i}$ is the length of the edge from the ideal vertex to the
  hyperideal vertex $i$, truncated at the horosphere and the dual plane.
\end{lemma}

Equations \eqref{eq:a_ij} and \eqref{eq:a_i} provide formulas for the edge
lengths of the truncated hyperbolic polyhedron in terms of the dihedral
angles $\alpha_{ij}$ and $\gamma_i$. Because the choice of the truncating
horosphere at the ideal vertex is arbitrary, the lengths $a_{i}$ are only
determined up to an additive constant. 

The hyperbolic lengths $a_{ij}$ and $a_{i}$ are related to the euclidean
lengths $l_{ij}$ and the radii $r_i$ (see Figure~\ref{fig:triangle}). The
radii $r_i$ are proportional to $e^{-a_i}$, \ie
\begin{equation}
  \label{eq:r_i_and_a_i}
  \frac{r_i}{r_j} = e^{-a_i + a_j},
\end{equation}
and
\begin{equation}
  \label{eq:l_ij_and_a_ij}
  l_{ij}^2 = r_1^2 + r_2^2 + 2 r_1 r_2 \cosh a_{12}.
\end{equation}
(These relations are obtained by straightforward calculation in the
Poin\-ca\-r\'e half-space model. The quantity $\cosh a_{12}$ is the
\emph{inversive distance} of two circles in the plane.) Together, Equations
\eqref{eq:a_ij}--\eqref{eq:l_ij_and_a_ij} provide formulas for the radii
$r_i$ and the euclidean edge lengths $l_{ij}$ in terms of the angles
$\alpha_{ij}$ and $\gamma_i$. They determine the $r_i$ and $l_{ij}$ up to a
common factor, in agreement with the fact that the angles determine the
decorated triangle of Figure~\ref{fig:triangle} up to similarity.

Now let $(\alpha,\gamma)\in\A(\T,\theta,\Xi)$ be a coherent angle
system. For each triangle $t\in T$ the angles
$(\alpha_{ij}^t,\alpha_{jk}^t,\alpha_{ki}^t,
\gamma_i^t,\gamma_j^t,\gamma_k^t)$ determine a decorated triangle up to
similarity. These fit together to form a hyperideal circle pattern iff they
can be scaled consistently; this means iff a radius $r_i$ can be assigned to
each vertex $i$ and a length $l_{ij}$ to each edge such that the relations
\eqref{eq:a_ij}--\eqref{eq:l_ij_and_a_ij} hold for each triangle.
Equivalently, the corresponding hyperbolic tetrahedra fit together iff the
horospheres at the infinite vertices can be chosen consistently; this means
iff an $a_i\in\R$ can be assigned to each vertex $i$ and an
$a_{ij}\in\R$ to each edge $ij$ such that the following holds:

If $t\in T$ and $i,j\in t$, then
\begin{equation}
  \label{eq:a_ij_condition}
  -2 \frac{\partial V}{\partial\alpha^t_{ij}} = a_{ij}
\end{equation}
and
\begin{equation}
  \label{eq:a_i_condition}
  -2 \Big(\frac{\partial V}{\partial\gamma^t_{i}}
  - \frac{\partial V}{\partial\gamma^t_{j}}\Big) = a_{i} - a_{j},
\end{equation}
where the derivatives are evaluated at 
$(\alpha^t_{ij},\alpha^t_{jk},\alpha^t_{ki},\gamma^t_i,\gamma^t_j,\gamma^t_k)$.
Now since 
\begin{equation*}
  \frac{\partial V}{\partial\alpha^t_{ij}}
  (\alpha^t_{ij},\alpha^t_{jk},\alpha^t_{ki},\gamma^t_i,\gamma^t_j,\gamma^t_k)
  = \frac{\partial F}{\partial\alpha^t_{ij}}(\alpha,\gamma)
\end{equation*}
and
\begin{equation*}
  \frac{\partial V}{\partial\gamma^t_{i}}
  (\alpha^t_{ij},\alpha^t_{jk},\alpha^t_{ki},\gamma^t_i,\gamma^t_j,\gamma^t_k)
  = \frac{\partial F}{\partial\gamma^t_{i}}(\alpha,\gamma),
\end{equation*}
(see Equation~\eqref{eq:F_sum}),
Equations \eqref{eq:a_ij_condition} and \eqref{eq:a_i_condition} are
equivalent to
\begin{equation}
  \label{eq:a_ij_condition_F}
  -2 \frac{\partial F}{\partial\alpha^t_{ij}}(\alpha,\gamma) = a_{ij}
\end{equation}
and
\begin{equation}
  \label{eq:a_i_condition_F}
  -2 \Big(\frac{\partial F}{\partial\gamma^t_{i}}(\alpha,\gamma)
  - \frac{\partial F}{\partial\gamma^t_{j}}(\alpha,\gamma)\Big) 
  = a_{i} - a_{j}.
\end{equation}
Clearly, the Equations \eqref{eq:a_ij_condition_F} for the $a_{ij}$ are
compatible iff the following condition holds:

\smallskip
(i) For each interior edge $ij$,
\begin{equation}
  \label{eq:compat_1}
  \Big(\frac{\partial}{\partial\alpha^t_{ij}}
  -\frac{\partial}{\partial\alpha^{t'}_{ji}}\Big)F(\alpha,\gamma)=0,
\end{equation}
where $t$ and $t'$ are the triangles adjacent with $ij$. 

Equations \eqref{eq:a_i_condition_F} for the $a_i$ are
compatible iff the condition holds: 

\smallskip
(ii) If $i_0 t_1 i_1 t_2 i_3 \ldots t_n i_n$ is any
finite sequence of alternatingly vertices and triangles that starts and ends
with the same vertex $i_0=i_n$, and that has the property that each $t_m$
contains preceding vertex $i_{m-1}$ and the following vertex $i_m$, then
\begin{equation}
  \label{eq:compat_2}
  \sum_{m=1}^n \Big(
  \frac{\partial}{\partial\gamma^{t_m}_{m}}-
  \frac{\partial}{\partial\gamma^{t_m}_{m-1}}
  \Big) F(\alpha,\gamma) = 0.
\end{equation}

These conditions (i) and (ii) are the non-linear compatibility conditions
that a coherent has to satisfy to represent a solution to the circle pattern
problem. It remains to show that they are also the conditions for a critical
point of $F$ under variations in $\A(\T,\theta,\Xi)$. This is achieved
by the following lemma, which concludes the proof of Lemma~\ref{lem:crit}.

\begin{lemma}
  The tangent space to $\A(\T,\theta,\Xi)$ is spanned by the tangent vectors
  \begin{equation}
    \label{eq:tangent_alpha}
    \frac{\partial}{\partial\alpha^t_{ij}}
    -\frac{\partial}{\partial\alpha^{t'}_{ji}}
  \end{equation}
  (one for each interior edge) and
  \begin{equation}
    \label{eq:tangent_gamma}
    \sum_{m=1}^n \Big(
    \frac{\partial}{\partial\gamma^{t_m}_{m}}-
    \frac{\partial}{\partial\gamma^{t_m}_{m-1}}
    \Big)
  \end{equation}
  (one for each cycle $i_0 t_1 i_1 t_2 i_3 \ldots t_n i_n$) that appear in
  conditions (i) and (ii) above.
\end{lemma}

\begin{proof}
  The space $\A(\T,\theta,\Xi)$ is defined by strict
  inequalities and Equations~\eqref{eq:gamma_triangle_sum},
  \eqref{eq:alpha_sum}, \eqref{eq:alpha_bdy},
  and~\eqref{eq:gamma_vertex_sum}. The equations for a tangent vector are
  therefore:

  \smallskip
  \begin{compactitem}
  \item[(a)] For each boundary edge $ij$:\quad $d\alpha^t_{ij}=0$. 
  \item[(b)] For each interior edge $ij=t\cap t'$:\quad 
    $d\alpha^t_{ij}+d\alpha^{t'}_{ji}=0$.
  \item[(c)] For each vertex $i$:\quad 
    $\displaystyle\sum_{t\ni i}d\gamma^t_i=0$.
  \item[(d)] For each triangle $t=ijk$:\quad
    $d\gamma^t_i+d\gamma^t_j+d\gamma^t_k=0$. 
  \end{compactitem}

  \smallskip

  Since each equation involves either alphas or gammas but not both,
  the tangent space is the direct sum of an $\alpha$-subspace and a
  $\gamma$-subspace. The $\alpha$-subspace is clearly spanned by the tangent
  vectors \eqref{eq:tangent_alpha}. To see that the $\gamma$-subspace is
  spanned by the vectors \eqref{eq:tangent_gamma}, let $G$ be
  the graph with vertex-set $\mathcal V=V\cup T$ and edge-set 
  \begin{equation*}
    \mathcal E=\big\{\{i,t\}\in\mathcal V \,\big|\, 
      i \in V, t\in T, i\in t\big\}.
  \end{equation*}
  The edges $\{i,t\}$ of $G$ are in one-to-one correspondence with the
  tangent vectors $\frac{\partial}{\partial\gamma^t_i}$. This gives rise to a
  linear isomorphism between the space of edge-chains of $G$ (over $\R$) and
  the space spanned by the vectors $\frac{\partial}{\partial\gamma^t_i}$.
  Equations (c) and (d) are then simply the equations for the cycle-space of
  $G$.
\end{proof}

\section{Proof of Lemma~\ref{lem:concave}}

\noindent
We are going to show that each of the terms
$V(\alpha_i^t,\alpha_j^t,\alpha_k^t,\gamma_i^t,\gamma_j^t,\gamma_k^t)$ in
Equation~\eqref{eq:F_sum} is concave. To this end, we split the
$15$-term sum in Equation~\eqref{eq:V} which defines $V$ into five parts; see
Equation~\eqref{eq:V_five_tetra}. Each part represents the volume of an ideal
tetrahedron (Theorem~\ref{thm:ideal_vol}), which is known to be concave
(Lemma~\ref{lem:ideal_vol_concave}).

\subsection{The volume of an ideal tetrahedron}
\label{sec:volume_ideal}

Consider a hyperbolic tetrahedron with all four vertices on the sphere at
infinity. For each vertex the sum of the interior dihedral angles at the
adjacent edges is $\pi$. This implies that the angles at opposite edges are
equal~\cite{milnor82} \cite{milnor94:_volume_in_hyp}. An ideal tetrahedron is
therefore determined by three angles in the set
\begin{equation}
  \label{eq:T0}
  \Delta_0=
  \big\{
    (\gamma_1,\gamma_2,\gamma_3)\in\R^3 \;\big|\; 
    \gamma_i>0,\,\gamma_1+\gamma_2+\gamma_3=\pi
  \big\}.
\end{equation}

\begin{theorem}[Milnor~\cite{milnor82} \cite{milnor94:_volume_in_hyp}]
  \label{thm:ideal_vol}
  The hyperbolic volume of an ideal tetrahedron with dihedral angles
  $(\gamma_1,\gamma_2,\gamma_3)\in\Delta_0$ is
  \begin{equation}
    \label{eq:Milnor}
    V_0(\gamma_1,\gamma_2,\gamma_3)=\ML(\gamma_1)+\ML(\gamma_2)+\ML(\gamma_3). 
  \end{equation}
\end{theorem}

\begin{lemma}[Rivin~\cite{rivin94}]
  \label{lem:ideal_vol_concave}
  The volume function $V_0$ is strictly concave on $\Delta_0$. 
\end{lemma}

\noindent The proof is straight forward. For the reader's convenience, we
repeat it here.

\begin{proof}[Proof of Lemma~\ref{lem:ideal_vol_concave}]
  Let 
  \begin{equation*}
    f(\alpha,\beta)=V_0(\alpha,\beta,\pi-\alpha-\beta)
    =\ML(\alpha)+\ML(\beta)-\ML(\alpha+\beta)
  \end{equation*}
and assume that $\alpha>0$, $\beta>0$, $\alpha+\beta<\pi$. Since
$\ML''(x)=-\cot x$, the Hessian matrix of $f$ is
\begin{equation*}
  \operatorname{Hess} f =
  \begin{pmatrix}
    -\cot\alpha+\cot(\alpha+\beta) & \cot(\alpha+\beta) \\
    \cot(\alpha+\beta) & -\cot\beta+\cot(\alpha+\beta)
  \end{pmatrix}.
\end{equation*}
A short calculation shows that the determinant of $\operatorname{Hess}(f)$ is
$1$. The matrix is therefore either positive definite or negative definite.
But since the cotangent is a strictly decreasing function on $(0,\pi)$, the
diagonal elements are negative. Hence $\operatorname{Hess}(f)$ is negative
definite and $f(\alpha,\beta)$ is strictly concave.
\end{proof}

\subsection{Five ideal tetrahedra}
\label{sec:five_ideal}

Equation~\eqref{eq:V} for the truncated volume
$V(\alpha_{12},\alpha_{23},\alpha_{31},\gamma_1,\gamma_2,\gamma_3)$ can be
rewritten as
\begin{equation*}
  2V(T)=\sum_{i=1}^3\big(
  \ML(\gamma_i)+\ML(\gamma_i')+\ML(\gamma_i'')+\ML(\mu_i)+\ML(\nu_i)
  \big),
\end{equation*}
where
\begin{equation*}
\begin{alignedat}{3}
  \gamma_1'&=\tfrac{\pi+\alpha_{31}-\alpha_{12}-\gamma_1}{2}, & \quad 
  \gamma_2'&=\tfrac{\pi+\alpha_{12}-\alpha_{23}-\gamma_2}{2}, & \quad 
  \gamma_3'&=\tfrac{\pi+\alpha_{23}-\alpha_{31}-\gamma_3}{2},\\
  \gamma_1''&=\tfrac{\pi-\alpha_{31}+\alpha_{12}-\gamma_1}{2}, & \quad 
  \gamma_2''&=\tfrac{\pi-\alpha_{12}+\alpha_{23}-\gamma_2}{2}, & \quad 
  \gamma_3''&=\tfrac{\pi-\alpha_{23}+\alpha_{31}-\gamma_3}{2}, \\
  \mu_1&=\tfrac{\pi+\alpha_{31}+\alpha_{12}-\gamma_1}{2}, & \quad
  \mu_2&=\tfrac{\pi+\alpha_{12}+\alpha_{23}-\gamma_2}{2}, & \quad
  \mu_3&=\tfrac{\pi+\alpha_{23}+\alpha_{31}-\gamma_3}{2}, \\
  \nu_1&=\tfrac{\pi-\alpha_{31}-\alpha_{12}-\gamma_1}{2}, & \quad
  \nu_2&=\tfrac{\pi-\alpha_{12}-\alpha_{23}-\gamma_2}{2}, & \quad
  \nu_3&=\tfrac{\pi-\alpha_{23}-\alpha_{31}-\gamma_3}{2}.
\end{alignedat}  
\end{equation*}
The following observation is both very simple and crucial for this proof (and
also for the proof of Lemma~\ref{lem:max} in
Section~\ref{sec:proof_lem_max}): If
\begin{equation*}
  (\alpha_{12},\alpha_{23},\alpha_{31},\gamma_1,\gamma_2,\gamma_3)\in\Delta,
\end{equation*}
then
\begin{equation*}
  (\gamma_1',\gamma_2',\gamma_3')\in\Delta_0,\qquad
  (\gamma_1'',\gamma_2'',\gamma_3'')\in\Delta_0
\end{equation*}
and also
\begin{equation*}
  (\gamma_1,\mu_1,\nu_1)\in\Delta_0,\qquad
  (\gamma_2,\mu_2,\nu_2)\in\Delta_0,\qquad
  (\gamma_3,\mu_3,\nu_3)\in\Delta_0.
\end{equation*}
Thus, $2V$ is the sum of the volumes of five ideal tetrahedra:
\begin{multline}
  \label{eq:V_five_tetra}
  2V(\alpha_{12},\alpha_{23},\alpha_{31},\gamma_1,\gamma_2,\gamma_3)=
  V_0(\gamma_1',\gamma_2',\gamma_3') 
  +V_0(\gamma_1'',\gamma_2'',\gamma_3'')\\
  +V_0(\gamma_1,\mu_1,\nu_1)+V_0(\gamma_2,\mu_2,\nu_2)
  +V_0(\gamma_3,\mu_3,\nu_3).
\end{multline}
Since each of the five terms is concave by Lemma~\ref{lem:ideal_vol_concave},
$2V$ is concave and so is $F$. This completes the proof of
Lemma~\ref{lem:concave}.

\begin{remark}
  We have no geometric explanation why $2V$ is the sum of the volumes of five
  ideal tetrahedra, although this may well be a consequence of Doyle \&
  Leibon's ``23040 symmetries of hyperbolic tetrahedra''~\cite{doyle03}.
  Equation~\eqref{eq:V_five_tetra} is not the only way to write $2V$ as a sum
  of five tetrahedra: For example, because
  $(\alpha_1,\alpha_2,\alpha_3,\gamma_1,\gamma_2,\gamma_3)\in\Delta$ also
  implies
  \begin{equation*}
    (\gamma_1,\gamma_2,\gamma_3)\in\Delta_0,\;
    (\gamma_1',\mu_2,\nu_3)\in\Delta_0,\;
    (\gamma_2',\mu_3,\nu_1)\in\Delta_0,\;
    (\gamma_3',\mu_1,\nu_2)\in\Delta_0,
  \end{equation*}
  one has as well
  \begin{multline*}
    2V(\alpha_{12},\alpha_{23},\alpha_{31},\gamma_1,\gamma_2,\gamma_3)=
    V_0(\gamma_1,\gamma_2,\gamma_3)
    +V_0(\gamma_1'',\gamma_2'',\gamma_3'')\\
    +V_0(\gamma_1',\mu_2,\nu_3)
    +V_0(\gamma_2',\mu_3,\nu_1)
    +V_0(\gamma_3',\mu_1,\nu_2).
  \end{multline*}
\end{remark}

\section{Proof of Lemma~\ref{lem:max}}
\label{sec:proof_lem_max}

\noindent
To prove Lemma~\ref{lem:max} we have to show the following:

\begin{claim*}
  Suppose $\A(\T,\theta,\Xi)\not=\emptyset$ and let
  $p\in\overline{\A(\T,\theta,\Xi)}\setminus\A(\T,\theta,\Xi)$. Then
  there is a $q\in\A(\T,\theta,\Xi)$ with $F(q)>F(p)$.
\end{claim*}

In the following Section, we will analyze the behavior of the volume function
$V$ as the dihedral angles approach the relative boundary of the domain. In
Section~\ref{sec:proof_claim} we will use this analysis to prove the Claim.

\subsection{Behavior of the volume function at the boundary of the domain}

We will first recollect Rivin's analysis~\cite{rivin94} of the behavior of
the volume $V_0(p)$ of an ideal tetrahedron (see
Section~\ref{sec:volume_ideal}) as $p$ approaches the relative boundary
$\overline{\Delta}_0\setminus\Delta_0$ of the domain $\Delta_0$. Then we will
use this for the corresponding analysis of the volume function
$V$ of a tetrahedron with one ideal an three hyperideal vertices. Here we
will again make essential use of the decomposition into five ideal tetrahedra
described in Section~\ref{sec:five_ideal}.

\subsubsection*{Boundary behavior of $V_0$}

First consider the set $\Delta_0$ of dihedral angles of ideal tetrahedra (see
Equation~\eqref{eq:T0} in Section~\ref{sec:volume_ideal}). Its closure
$\overline{\Delta}_0$ is the $2$-simplex in $\R^3$ that is spanned by the
points $(\pi,0,0)$, $(0,\pi,0)$, $(0,0,\pi)$. The points in the relative
boundary $\overline{\Delta}_0\setminus\Delta_0$ correspond to ideal
tetrahedra that have degenerated to planar figures.

\begin{definition}
  We call a point
  $(\gamma_1,\gamma_2,\gamma_3)\in\overline{\Delta}_0\setminus\Delta_0$
  \emph{mildly degenerate} iff $(\gamma_1,\gamma_2,\gamma_3)$ is some
  permutation of $(0,\beta,\pi-\beta)$ with $0<\beta<\pi$, \ie{} iff $p$ is
  contained in an open side of the boundary triangle. We call it \emph{badly
    degenerate} iff $(\gamma_1,\gamma_2,\gamma_3)$ is some permutation of
  $(0,0,\pi)$, \ie{} iff it is a vertex of the boundary triangle.
\end{definition}

It is easy to see that $V_0$ vanishes on
$\overline{\Delta}_0\setminus\Delta_0$. But the speed with which $V(q)$ tends
to $0$ as $q\in\Delta_0$ approaches the boundary is different depending on
whether $q$ approaches a mildly or a badly degenerate point: 

\begin{lemma}
  \label{lem:boundary_derivative}
  Let $p\in\overline{\Delta}_0\setminus\Delta_0$, $q\in\Delta_0$. If $p$ is
  mildly degenerate, then
  \begin{equation}
    \label{eq:infinite_derivative}
    \lim_{t\searrow 0}\frac{d}{dt}V_0\big((1-t)p+t q\big) = +\infty.
  \end{equation}
  If $p$ is badly degenerate then the $t$-derivative
  \begin{equation}
    \label{eq:finite_derivative}
    \frac{d}{dt}V_0\big((1-t)p+t q\big)
  \end{equation}
  has a finite positive limit for $t\searrow 0$.
\end{lemma}
\begin{proof}
  The claim for mildly degenerate $p$ follows directly from the fact that
  $\ML(x)$ is smooth except at integer multiples of $\pi$, where the
  derivative $\ML'(x)=-\log|2\sin x|$ tends to $+\infty$. To prove
  the claim for badly degenerate $p$, let us assume without loss of
  generality that $p=(0,0,\pi)$. Then
  \begin{equation*}
    (1-t)p+t q = (ta,tb,\pi-t(a+b))
  \end{equation*}
  for some $a,b>0$, and
  \begin{equation*}
    V_0((1-t)p+t q)=\ML(ta)+\ML(tb)-\ML\big(t(a+b)\big).
  \end{equation*}
  Hence for the $t$-derivative we obtain
  \begin{equation*}
    \frac{d}{dt}V_0\big((1-t)p+t q\big)=
    \log\bigg|
    \frac{\sin^{a+b}\big(t(a+b)\big)}{\sin^a(ta)\sin^b(tb)}
    \bigg|
    \xrightarrow{t\searrow 0}\log\frac{(a+b)^{a+b}}{a^a b^b}>0.
  \end{equation*}
\end{proof}

\subsubsection*{Boundary behavior of $V$} Now consider the set $\Delta$ of
dihedral angles of tetrahedra with one ideal and three hyperideal vertices
(see Section~\ref{sec:local_geom}). Its closure $\overline{\Delta}$ is the
set of points
$(\alpha_{12},\alpha_{23},\alpha_{31},\gamma_1,\gamma_2,\gamma_3)\in\R^6$
that satisfy Equation~\eqref{eq:gamma_triangle_sum} and the non-strict
versions of the inequalities~\eqref{eq:alpha_gamma_positive} and
\eqref{eq:alpha_gamma_inequalities}. A point is in the relative boundary
$\overline{\Delta}\setminus\Delta$ if it satisfies with equality at least one
of these the non-strict inequalities. In Section~\ref{sec:five_ideal} we
described an affine map $\Delta\rightarrow(\Delta_0)^5$ associating five
ideal tetrahedra with each point in $\Delta$. This extends to a map
$\overline{\Delta}\rightarrow(\overline{\Delta}_0)^5$. We classify the
degenerate points $p\in\overline{\Delta}\setminus\Delta$ according to whether
any of the five corresponding ideal tetrahedra degenerate and the way in
which they do:

\begin{definition}
  Let $p\in\overline{\Delta}\setminus\Delta$. We say that $p$ is \emph{mildly
    degenerate} iff at least one of the five corresponding ideal tetrahedra
  is mildly degenerate. We say that $p$ is \emph{badly degenerate} iff at
  least one of the five ideal tetrahedra is badly degenerate but none are
  mildly degenerate. We say that $p$ is \emph{$\alpha$-degenerate} iff all
  five corresponding ideal tetrahedra are non-degenerate.
\end{definition}

Clearly, every $p\in\overline{\Delta}\setminus\Delta$ is either mildly
degenerate, badly degenerate or $\alpha$-degenerate. By
Lemma~\ref{lem:badly_degenerate}, there is only one badly degenerate $p$ up to
a permutation of the indices. The reason for calling the last type
``$\alpha$-degenerate'' will be made clear by
Lemma~\ref{lem:alpha_degenerate}.  The next lemma follows immediately from
Lemma~\ref{lem:boundary_derivative}.

\begin{lemma}
  \label{lem:boundary_derivative_of_V}
  Let $p\in\overline{\Delta}\setminus\Delta$, $q\in\Delta$. If $p$ is
  mildly degenerate, then
  \begin{equation}
    \lim_{t\searrow 0}\frac{d}{dt}V\big((1-t)p+t q\big) = +\infty.
  \end{equation}
  If $p$ is badly degenerate or $\alpha$-degenerate then the $t$-derivative
  \begin{equation}
    \frac{d}{dt}V\big((1-t)p+t q\big)
  \end{equation}
  has a finite limit for $t\searrow 0$.
\end{lemma}

\begin{lemma}
  \label{lem:badly_degenerate}
  Suppose $(\alpha_{12},\alpha_{23},\alpha_{31},\gamma_1,\gamma_2,\gamma_3)
  \in\overline{\Delta}\setminus\Delta$ is badly degenerate. Then
  \begin{equation*}
    \gamma_i=\alpha_{jk}=\pi,\qquad\gamma_j=\gamma_k=\alpha_{ij}=\alpha_{ki}=0
  \end{equation*}
  for some permutation $(i,j,k)$ of $(1,2,3)$. (This implies that \emph{all
    five} corresponding ideal tetrahedra are badly degenerate.)
\end{lemma}
\begin{proof}
  We consider separately all essentially different cases of one of the five
  ideal tetrahedra badly degenerating in some way, where we consider cases as
  essentially different if they do not only differ by a permutation of the
  indices. First, note that by summing (the non-strict versions of) the
  inequalities~\eqref{eq:alpha_gamma_inequalities} and using
  Equation~\eqref{eq:gamma_triangle_sum}, one obtains
  \begin{equation*}
    \alpha_{12}+\alpha_{23}+\alpha_{31}\leq\pi.
  \end{equation*}
  
  \emph{Case 1:} $(\gamma_1',\gamma_2',\gamma_3')=(0,0,\pi)$. Since
  \begin{equation*}
    \pi=2\gamma_3'-\pi=\alpha_{23}-\alpha_{31}-\gamma_3,
  \end{equation*}
  we have $\alpha_{23}=\pi$ and hence
  $\alpha_{12}=\alpha_{31}=\gamma_2=\gamma_3=0$ and $\gamma_1=\pi$.
  
  \emph{Case 2:} $(\gamma_1,\mu_1,\nu_1)=(\pi,0,0)$. First, $\gamma_1=\pi$
  implies $\gamma_2=\gamma_3=\alpha_{12}=\alpha_{31}=0$. Then
  \begin{equation*}
    (\gamma_1',\gamma_2',\gamma_3')= 
    (0,\frac{\pi-\alpha_{23}}{2},\frac{\pi+\alpha_{23}}{2}),
  \end{equation*}
  and since by assumption this is not mildly degenerate, we have
  $\alpha_{23}=\pi$. 
  
  \emph{Case 3:} $(\gamma_1,\mu_1,\nu_1)=(0,\pi,0)$. 
  First, $\gamma_1=0$ implies $\gamma_2+\gamma_3+\pi$; $\mu_1=\pi$ implies
  $\alpha_{31}+\alpha_{12}=\pi$ and hence $\alpha_{23}=0$. Now the non-strict
  versions of the inequalities~\eqref{eq:alpha_gamma_inequalities} imply
  $\alpha_{12}+\gamma_2=\pi$ and $\alpha_{31}+\gamma_3=\pi$. Hence
  $\nu_2=\nu_3=0$. Since by assumption, neither $(\gamma_2,\mu_2,\nu_2)$ nor
  $(\gamma_3,\mu_3,\nu_3)$ are mildly degenerate, either
  $\gamma_2=\alpha_{31}=\pi$ and $\gamma_3=\alpha_{12}=0$ or
  $\gamma_2=\alpha_{31}=0$ and $\gamma_3=\alpha_{12}=\pi$.

  \emph{(Non-)Case 4:} $(\gamma_1,\mu_1,\nu_1)=(0,0,\pi)$. This cannot happen
  because $\nu_1\leq\frac{\pi}{2}$.
\end{proof}

\begin{lemma}
  \label{lem:alpha_degenerate}
  Suppose $p=(\alpha_{12},\alpha_{23},\alpha_{31},\gamma_1,\gamma_2,\gamma_3)
  \in\overline{\Delta}\setminus\Delta$ is $\alpha$-degenerate. Then all
  $\gamma_i>0$ and the strict inequalities~\eqref{eq:alpha_gamma_inequalities}
  are satisfied, but some $\alpha_{ij}$ vanish.
\end{lemma}
\begin{proof}
  First, $\gamma_i>0$ and $\gamma_i+\alpha_{ij}+\alpha_{ki}<\pi$ because
  $\gamma_i=0$ or $\gamma_i+\alpha_{ij}+\alpha_{ki}=\pi$ would imply that
  $(\gamma_i,\mu_i,\nu_i)$ is degenerate. But then some $\alpha_{ij}$ must
  vanish because otherwise $p\in\Delta$.
\end{proof}

\subsection{Proof of the Claim}
\label{sec:proof_claim}

Suppose $\A(\T,\theta,\Xi)\not=\emptyset$ and let
$p=(\alpha,\gamma)\in\overline{\A(\T,\theta,\Xi)}\setminus\A(\T,\theta,\Xi)$.
We distinguish several cases.

\bigskip\noindent
\textit{Case 1. There is at least one triangle $t\in T$ such that
  $(\alpha^t,\gamma^t)$ is mildly degenerate.}
Let $\tilde p=(\tilde\alpha,\tilde\gamma)\in\A(\T,\theta,\Xi)$ be any
coherent angle system. Then Lemma~\ref{lem:boundary_derivative_of_V} implies
\begin{equation*}
  \lim_{t\searrow 0}\frac{d}{dt}F\big((1-t)p+t\tilde p\big) = +\infty.
\end{equation*}
It follows (by the mean value theorem) that if $\varepsilon>0$ is small
enough, then $F(q)>F(p)$ for $q=(1-\varepsilon)p+\varepsilon\tilde p$. This
completes the proof of the Claim under the assumption of Case~1.

\bigskip\noindent
\textit{Case 2. There are no $t\in T$ with $(\alpha^t,\gamma^t)$
  mildly degenerate or $\alpha$-degenerate.} 
This means all degenerate $(\alpha^t,\gamma^t)$ are badly degenerate. We will
construct a $\tilde
p\in\overline{\A(\T,\theta,\Xi)}\setminus\A(\T,\theta,\Xi)$ which satisfies
$F(p)=F(\tilde p)$ and the conditions of Case~1. Since the Claim was proven
for Case~1, it holds in Case~2 also.

\begin{figure}
  \centering \input{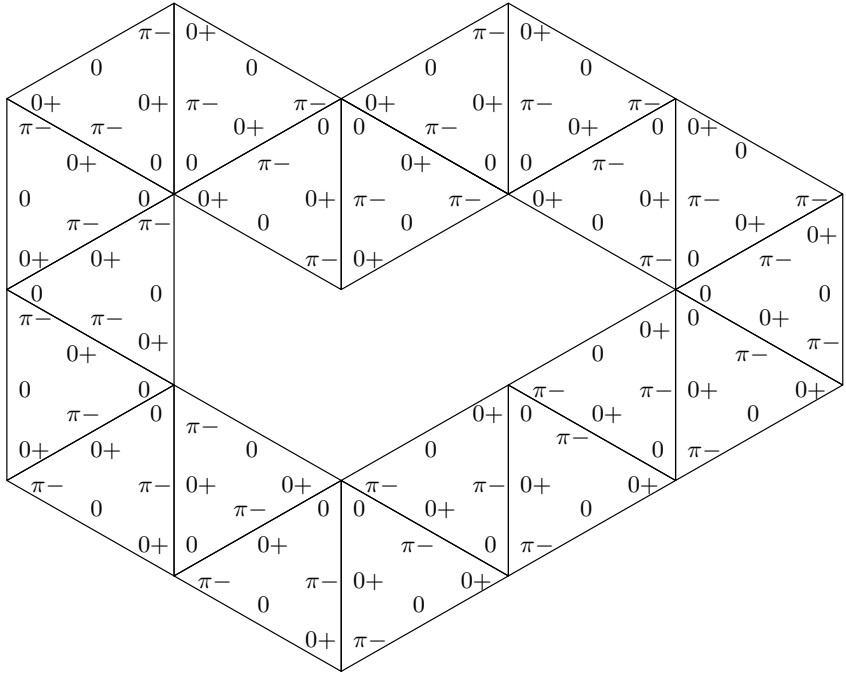}
  \caption{How to modify the angles along a loop of badly degenerate
    triangles. We write $\pi-$ and $0+$ as shorthand for ``$\pi$ is replaced
    by $\pi-\varepsilon$'' and for ``$0$ is replaced by $\varepsilon$'',
    respectively. The changes in $\alpha$-angles are indicated along the
    edges, the changes in $\gamma$-angles are indicated in the corners of
    each triangle.}
  \label{fig:deform_badly}
\end{figure}

Suppose $t_1\in T$ with $(\alpha^t,\gamma^t)$ badly degenerate. Let $ij$ be
the edge of $t_1$ with $\alpha^{t_1}_{ij}=\pi$.
Equations~\eqref{eq:alpha_sum} and~\eqref{eq:alpha_bdy} imply that
$\theta_{ij}=0$ and therefore that edge $ij$ is not a boundary edge. Let
$t_2\in T$ be the triangle neighboring $t_1$ across edge $ij$. Again by
Equation~\eqref{eq:alpha_sum}, $\alpha^{t_2}_{ji}=0$, so $t_2$ is also badly
degenerate. By repeating this argument we construct a sequence
$t_1,t_2,t_3,\ldots$ of badly degenerate triangles, each one adjacent to the
next. Since there are only finitely many triangles, this sequence must
eventually loop back on itself.  Let us reindex the triangles such that
$t_1,t_2,\ldots,t_n=t_1$ is such a loop of badly degenerate triangles. Define
an angle system $\tilde p=(\tilde\alpha,\tilde\gamma)$ as follows (see
Figure~\ref{fig:deform_badly}). 
If $t\in T$ is not contained in the loop of
triangles, then let $(\tilde\alpha^t,\tilde\gamma^t)=(\alpha^t,\gamma^t)$. If
$t=t_m$ is contained in the loop, let $ij=t_m\cap t_{m+1}$ and
$jk=t_{m-1}\cap t_m$, hence
\begin{equation*}
  \alpha^t_{ij}=\gamma^t_k=\pi,\qquad
  \alpha^t_{jk}=\alpha^t_{ki}=\gamma_i=\gamma_j=0.
\end{equation*}
Let 
\begin{align*}
  \tilde\alpha^t_{ij}&=\tilde\gamma^t_k=\pi-\varepsilon,\\
  \tilde\alpha^t_{jk}&=\tilde\gamma_i=\varepsilon,\\
  \tilde\alpha^t_{ki}&=\tilde\gamma_j=0,
\end{align*}
with $\varepsilon\in(0,\pi)$ arbitrary but the same for all triangles in the
loop.  We claim that $\tilde
p\in\overline{\A(\T,\theta,\Xi)}\setminus\A(\T,\theta,\Xi)$. Indeed, the
triangles in the loop are still in $\overline\Delta\setminus\Delta$, but now
they are mildly degenerate instead of badly degenerate. Further, the sum of
$\alpha$-angles at edges has obviously not changed, so
Equation~\eqref{eq:alpha_sum} still holds. Finally, to see that the sum of
$\gamma$-angles around each vertex has not changed, note that at a vertex $i$
this angle sum is increased by $\varepsilon$ for each time that the loop of
triangles enters the star of $i$ and is decreased by $\varepsilon$ each time
it leaves the star. So Equation~\eqref{eq:gamma_vertex_sum} still holds. Also
$F(p)=F(\tilde p)$, because the truncated volumes for the triangles in the
loop are $0$ before and after the deformation. Hence we have reduced Case~2 to
Case~1.

\bigskip\noindent
\textit{Case 3. There are no $t\in T$ with
  $(\alpha^t,\gamma^t)$ mildly degenerate, but there is at least one $t\in T$
  such that $(\alpha^t,\gamma^t)$ is $\alpha$-degenerate.}
Suppose $t=ijk\in T$ is $\alpha$-degenerate and $\alpha^t_{ij}=0$. Let
$t'=jil$ be the triangle on the other side of edge $ij$. Note that
$\alpha^{t'}_{ji}$ cannot vanish, because this would imply
$\theta_{ij}=\pi-\alpha^t_{ij}-\alpha^{t'}_{ji}=\pi$. We distinguish two
sub-cases.

\medskip\noindent
\textit{Case 3(a). $t'$ is not badly degenerate.} 
Then $0<\alpha^{t'}_{ji}<\pi$, and we can change the angle system $p$ to the
angle system $\tilde p$ by setting $\tilde\alpha^t_{ij}=\varepsilon$ and
$\tilde\alpha^{t'}_{ji}=\alpha^{t'}_{ji}-\varepsilon$ for some small
$\varepsilon>0$, and keeping all other angles the same.  We claim that
$F(\tilde p)>F(p)$ if $\varepsilon$ is small enough. Indeed,
\begin{equation*}
  \frac{\partial}{\partial\alpha^t_{ij}}V(\alpha^t,\gamma^t)=0,
\end{equation*}
because $V(\alpha^t,\gamma^t)$ is even in the $\alpha$-variables (see
Equation~\eqref{eq:V}), and
\begin{equation*}
  \frac{\partial}{\partial\alpha^{t'}_{ji}}V(\alpha^{t'},\gamma^{t'})=
  -\frac{1}{2}a^{t'}_{ji}<0,
\end{equation*}
where $a^{t'}_{ji}$ is the length of the corresponding truncated edge
(Lemma~\ref{lem:V_derivative}), which is positive because $a_{ji}=0$
only if $\alpha_{ji}=0$ (see Equation~\eqref{eq:l_ij_and_a_ij} and
Figure~\ref{fig:triangle}). So provided that $\varepsilon>0$ is small
enough, $V(\alpha^t,\gamma^t)+V(\alpha^{t'},\gamma^{t'})$ increases if
$\alpha^t_{ij}$ increases by $\varepsilon$ and $\alpha^{t'}_{ji}$
decreases by $\varepsilon$. We have thus constructed a $\tilde
p\in\overline{\A(\T,\theta,\Xi)}$ without mildly degenerated triangles
such that $F(\tilde p)>F(p)$ and the number of vanishing
$\alpha$-angles has decreased by one. The Claim follows by induction
on the number of vanishing $\alpha$-angles.

\begin{figure}
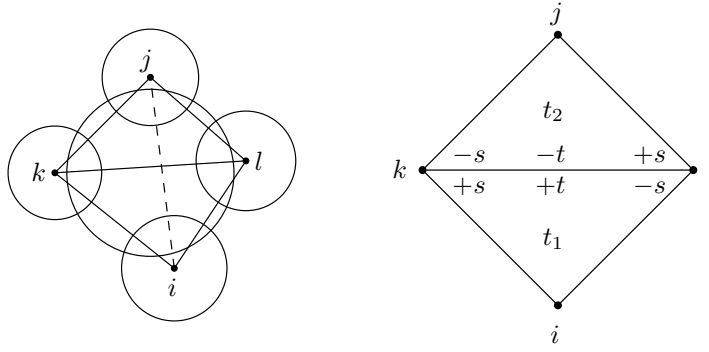

  \centering
  \raisebox{0.4cm}{\input{figs/second_flip.tex}}
  \hspace{1cm}
  \input{figs/deform_tilde_q.tex}
  \caption{\emph{Left:} If neighboring triangles share the same orthogonally
    intersecting circle, we can perform an edge flip. \emph{Right:} How
    $\tilde q(s,t)$ differs from $\tilde q$.}
  \label{fig:second_flip+deform_tilde_q}
\end{figure}

\medskip\noindent
\textit{Case 3(b). $t'$ is badly degenerate} 
We will construct a $q\in\A(\T,\theta,\Xi)$ such that $F(q)>F(p)$, but the
construction is less straightforward than in the other cases. Note that in
this case, $\alpha^{t'}_{ji}=\pi$ and hence $\theta_{ij}=0$. This means that
the solution of the circle pattern problem (if it exists) is such that the
same face-circle corresponds to both triangles $t$ and $t'$; see
Figure~\ref{fig:second_flip+deform_tilde_q} (left).
Equivalently, the tetrahedra corresponding to these triangles fit together
to form a pyramid over a quadrilateral base. Thus, one could pose an
equivalent circle pattern problem using the triangulation $\tilde\T$ that is
obtained by flipping the edge $ij$. This is the basic motivating idea behind
the construction we will now describe.

\begin{figure}
  \centering
  \input{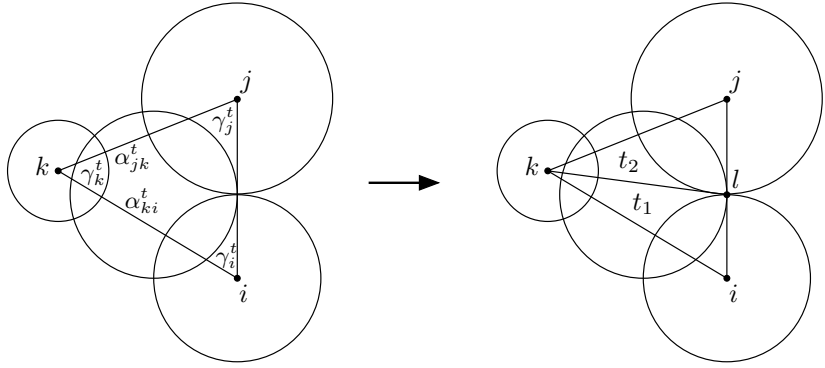}
  \caption{The $\alpha$-degenerate tetrahedron $t=(ijk)$ with
  $\alpha^t_{ij}=0$ is split into two triangles $t_1=(ilk)$, $t_2=(jkl)$.}
  \label{fig:deform_alpha}
\end{figure}

First, we split the $\alpha$-degenerate tetrahedron $(\alpha^t,\gamma^t)$
into two tetrahedra $(\alpha^{t_1},\gamma^{t_1})$ and
$(\alpha^{t_2},\gamma^{t_2})$ as shown in Figure~\ref{fig:deform_alpha}.
This introduces new angles $\gamma^{t_1}_k$, $\gamma^{t_2}_k$,
$\gamma^{t_1}_l$, $\gamma^{t_2}_l$, $\alpha^{t_1}_{lk}$, $\alpha^{t_2}_{kl}$,
which are uniquely determined by the old angles $(\alpha^t, \gamma^t)$. They
clearly satisfy
\begin{equation*}
  \gamma^{t_1}_k+\gamma^{t_2}_k=\gamma^t_k, \quad
  \gamma^{t_1}_l+\gamma^{t_2}_l=\pi, \quad
  \alpha^{t_1}_{lk}+\alpha^{t_2}_{kl}=\pi,
\end{equation*}
and also 
\begin{equation*}
  \gamma^{t_1}_l+\alpha^{t_1}_{lk}=\pi,\qquad 
  \gamma^{t_2}_l+\alpha^{t_2}_{kl}=\pi.
\end{equation*}
The remaining angles in $(\alpha^{t_1},\gamma^{t_1})$ and
$(\alpha^{t_2},\gamma^{t_2})$ are equal to the corresponding angles in
$(\alpha^t, \gamma^t)$: 
\begin{equation*}
  \gamma^{t_1}_i=\gamma^t_i,\quad 
  \gamma^{t_2}_j=\gamma^t_j,\quad
  \alpha^{t_1}_{ki}=\alpha^t_{ki},\quad
  \alpha^{t_2}_{jk}=\alpha^t_{jk},\quad
  \alpha^{t_1}_{il}=\alpha^{t_2}_{lj}=\alpha^t_{ij}=0.
\end{equation*}
The volumes satisfy
\begin{equation}
  \label{eq:split_volume_sum}
  V(\alpha^t, \gamma^t)=V(\alpha^{t_1},\gamma^{t_1})+
  V(\alpha^{t_2},\gamma^{t_2}).
\end{equation}
The tetrahedra $(\alpha^{t_1},\gamma^{t_1})$ and
$(\alpha^{t_2},\gamma^{t_2})$ are mildly degenerate, because they are not
badly degenerate but
\begin{equation*}
  \gamma^{t_1}_l+\alpha^{t_1}_{lk}+\alpha^{t_1}_{il}=\pi,\quad
  \gamma^{t_2}_l+\alpha^{t_2}_{lj}+\alpha^{t_2}_{kl}=\pi.
\end{equation*}
Now let $\tilde\T$ be the triangulation obtained from $\T$ by flipping the
edge $ij$, thus replacing it with an edge $kl$ and replacing the triangles
$t$, $t'$ with triangles $t_1$, $t_2$. (This edge flip can be performed even
if the triangulation $\T$ is not regular. The only obstruction for an edge to
be flippable is that it is adjacent to the same triangle on either side. But
this is not the case here, because $(\alpha^t,\gamma^t)$ is
$\alpha$-degenerate, whereas $(\alpha^{t'},\gamma^{t'})$ is badly
degenerate.) Let $\tilde p$ be the angle system with
$(\alpha^{t_1},\gamma^{t_1})$, $(\alpha^{t_2},\gamma^{t_2})$ as described
above and all other angles the same as in $p$. Then
\begin{equation*}
  \tilde p\in\overline{\A(\tilde\T,\tilde\theta,\Xi)},
\end{equation*}
where $\tilde\theta_{mn}=\theta_{mn}$ for all edges $mn$ of $\tilde\T$ except
for $kl$ and $\tilde\theta_{kl}=\theta_{ij}=0$. Because of
Equation~\eqref{eq:split_volume_sum} and because the volume of the badly
degenerate tetrahedron $(\alpha^{t'},\gamma^{t'})$ vanishes, we have
$F_{\T}(p)=F_{\tilde\T}(\tilde p)$. Now
$\A(\tilde\T,\tilde\theta,\Xi)\not=\emptyset$ because from every coherent
angle system in $\A(\T,\theta,\Xi)$ one can easily construct a coherent angle
system in $\A(\tilde\T,\tilde\theta,\Xi)$. Hence the reasoning of Case~1
above applies and there exists a
\begin{equation*}
  \tilde q=(\tilde\alpha,\tilde\gamma)\in\A(\tilde\T,\tilde\theta,\Xi)
  \quad\text{with}\quad
  F_{\tilde\T}(\tilde q)> F_{\tilde\T}(\tilde p)=F_{\T}(p).
\end{equation*}
It remains to construct a $q\in\A(\T,\theta,\Xi)$ with $F_{\T}(q)\geq
F_{\tilde\T}(\tilde q)$. If the tetrahedra
$(\tilde\alpha^{t_1},\tilde\gamma^{t_1})$ and
$(\tilde\alpha^{t_2},\tilde\gamma^{t_2})$ fit together, this could be
achieved by performing another edge flip as in
Figure~\ref{fig:second_flip+deform_tilde_q} (left). However, they will in
general not fit together. We will therefore deform $\tilde q$ to obtain a
$\tilde{\tilde q}\in\A(\tilde\T,\tilde\theta,\Xi)$ such that the triangles
fit together and $F_{\tilde\T}(\tilde{\tilde q})\geq F_{\tilde\T}(\tilde q)$.
To this end, let $\tilde q(s,t)$ be the angle system with
\begin{equation*}
  \begin{alignedat}{2}
    \tilde\gamma^{t_1}_k(s,t)&=\tilde\gamma^{t_1}_k+s,&\qquad
    \tilde\gamma^{t_1}_l(s,t)&=\tilde\gamma^{t_1}_l-s,\\
    \tilde\gamma^{t_2}_k(s,t)&=\tilde\gamma^{t_2}_k-s,&\qquad
    \tilde\gamma^{t_2}_l(s,t)&=\tilde\gamma^{t_2}_l+s,\\
    \tilde\alpha^{t_1}_{lk}(s,t)&=\tilde\alpha^{t_1}_{lk}+t,&\qquad
    \tilde\alpha^{t_2}_{kl}(s,t)&=\tilde\alpha^{t_2}_{kl}-t,  
  \end{alignedat}
\end{equation*}
and all other angles the same as in $\tilde q$; see
Figure~\ref{fig:second_flip+deform_tilde_q} (right).  Let
\begin{equation*}
  U=\{(s,t)\in\R^2\,|\,\tilde q(s,t)\in\A(\tilde\T,\tilde\theta,\Xi)\}.
\end{equation*}
This is a bounded open subset of $\R^2$. Let $f(s,t)=F_{\tilde\T}(\tilde
q(s,t))$. For $(s_0, t_0)\in U$, the two tetrahedra
$(\tilde\alpha^{t_1}(s_0,t_0),\tilde\gamma^{t_1}(s_0,t_0))$ and
$(\tilde\alpha^{t_2}(s_0,t_0),\discretionary{}{}{}\tilde\gamma^{t_2}(s_0,t_0))$ fit together iff
$(s_0, t_0)$ is a critical point of $f(s,t)$. To see this, apply the same
reasoning that was used to prove Lemma~\ref{lem:crit} in
Section~\ref{sec:proof_lem_crit}. Also, $f(s,t)$ is strictly concave on $U$,
because it is a restriction of the strictly concave function $F_{\tilde\T}$
(Lemma~\ref{lem:concave}) to an affine subspace. Finally, the maximum of
$f(s,t)$ on the compact set $\overline U$ cannot be attained on the boundary
$\partial U$. To see this, apply the same reasoning that was used in Case~1
and in Case~3(a) above. (The tetrahedra
$(\tilde\alpha^{t_1}(s,t),\tilde\gamma^{t_1}(s,t))$ and
$(\tilde\alpha^{t_2}(s,t),\tilde\gamma^{t_2}(s,t))$ cannot degenerate badly
for $(s,t)\in\overline U$.) Hence the restriction of $f(s,t)$ to $\overline
U$ attains its maximum at some $(s_m,t_m)\in U$, and we have found
$\tilde{\tilde q}=\tilde q(s_m,t_m)$. Finally we obtain
$q\in\A(\T,\theta,\Xi)$ with $F_{\T}(q)=F_{\tilde\T}(\tilde{\tilde
  q})>F_{\T}(p)$ from $\tilde{\tilde q}$ by flipping the edge $kl$.

Since we have thus dealt with the last Case, this completes the proof of the
Claim, and hence the proof of Lemma~\ref{lem:max}.

\section{Volume computations. Proof of Lemma~\ref{lem:volume_V}}
\label{sec:volumes}

\noindent
In this section we derive Equation~\eqref{eq:V} for the volume of a
tetrahedron with one ideal and three hyperideal vertices, that is, we prove
Lemma~\ref{lem:volume_V} of Section~\ref{sec:proof_lem_crit}.
Vinberg~\cite{vinberg91} \cite{vinberg93} derived a formula---in fact, the
same formula---for the volume of a tetrahedron with one ideal and three
\emph{finite}\/ vertices. Here we follow a very similar path. We subdivide a
tetrahedron with one ideal and three hyperideal vertices into three special
pyramids (Section~\ref{sec:vol_1_ideal_3_hyper}). A volume formula for these
is derived in Section~\ref{sec:spec_pyramid} with the help of yet other volume
formulas that we present in the following sections.

\subsection{The volume of a birectangular tetrahedron with two ideal vertices}

\begin{figure}
  \centering
  \input{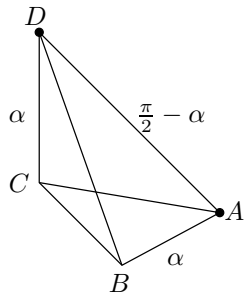}
  \caption{The birectangular tetrahedron $P_1$ with two ideal vertices $A$
    and $D$. The dihedral angles at the unlabeled edges are~$\frac{\pi}{2}$.}
  \label{fig:birectangular}
\end{figure}

A tetrahedron with vertices $ABCD$ is called \emph{birectangular} or an
\emph{orthoscheme} if the edge $AB$ is perpendicular on the side $BCD$ and
the edge $CD$ is perpendicular on the side $ABC$. Then the dihedral angles at
edges $AC$, $BD$, and $BC$ are $\frac{\pi}{2}$. A formula for the volume of a
birectangular hyperbolic tetrahedron as a function of the remaining three
dihedral angles was already derived by Lobachevsky~\cite{lobachevsky04}, see
also Coxeter~\cite{coxeter35}. We are only interested in the case of a
birectangular tetrahedron $P_1$ whose vertices $A$ and $D$ are ideal (see
Figure~\ref{fig:birectangular}). 
Because the dihedral angles sum to $\pi$ at the ideal vertices it follows
that the angles at $AB$ and $CD$ are equal, say, to $\alpha$, and the angle
at $AD$ is $\frac{\pi}{2}-\alpha$. Milnor derived the particularly simple
volume formula
\begin{equation}
  \label{eq:vol_P1}
  \Vol(P_1)=\frac{1}{2}\ML(\alpha)
\end{equation}
by direct integration~\cite{milnor82} \cite{milnor94:_volume_in_hyp}.

\subsection{The volume of an ideal prism}

Let $P_2$ be a triangular prism with all vertices at infinity. Such a prism is
always symmetric with respect to a planar reflection that interchanges the
triangular faces. (This is so because any three points on $S^2$, the sphere
at infinity, can be mapped to any other three points on $S^2$ by an
orientation reversing M\"obius transformation of $S^2$, and such a
transformation is the restriction of a hyperbolic reflection.)
Formula~\eqref{eq:vol_prism} below for the volume of $P$ is derived
in~\cite{leibon02}. For the reader's convenience we reproduce the argument.
Let the interior dihedral angles be $\alpha$, $\beta$, $\gamma$,
$\alpha'$, $\beta'$, $\gamma'$, as shown in
\begin{figure}[tbp]
  \centering
  \input{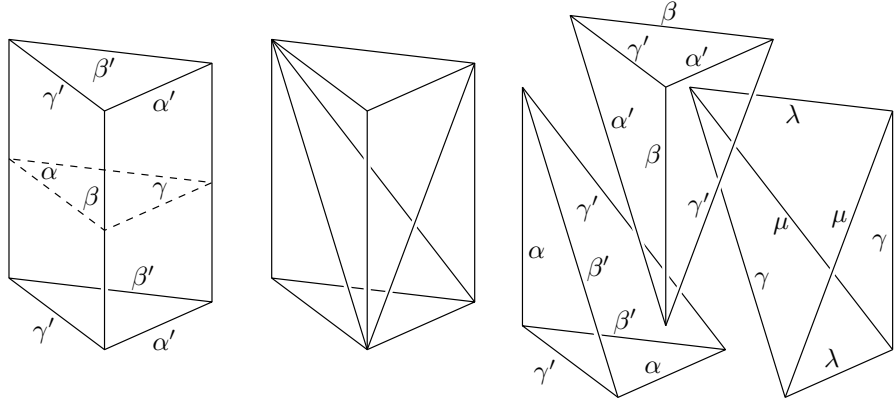}
  \caption{The ideal prism $P_2$. To calculate its volume, we subdivide it
    into three ideal tetrahedra.}
  \label{fig:prism}
\end{figure}
Figure~\ref{fig:prism} (left). The symmetry plane intersects the side faces
of the prism orthogonally in the dashed triangle. Hence there exists an ideal
prism with dihedral angles $\alpha$, $\beta$, $\gamma$, iff there is a
hyperbolic triangle with these angles, \ie{} iff
\begin{equation*}
  \alpha + \beta + \gamma < \pi.
\end{equation*}
The dihedral angles sum to $\pi$ at each ideal vertex, hence
\begin{equation*}
  \gamma' = \frac{\pi-\alpha-\beta+\gamma}{2},\quad
  \alpha' = \frac{\pi+\alpha-\beta-\gamma}{2},\quad
  \beta' = \frac{\pi-\alpha+\beta-\gamma}{2}\,.
\end{equation*}

\begin{lemma}[Leibon \cite{leibon02}]
  The volume of the prism $P_2$ is
  \begin{multline}
    \label{eq:vol_prism}
    \Vol(P_2) = \ML(\alpha) + \ML(\beta) + \ML(\gamma)\\
    + \ML\big(\frac{\pi+\alpha-\beta-\gamma}{2}\big)
    + \ML\big(\frac{\pi-\alpha+\beta-\gamma}{2}\big)\\
    + \ML\big(\frac{\pi-\alpha-\beta+\gamma}{2}\big)
    + \ML\big(\frac{\pi-\alpha-\beta-\gamma}{2}\big).
  \end{multline}  
\end{lemma}

\begin{proof}
  Subdivide the the prism $P_2$ into three ideal tetrahedra as shown in
  Figure~\ref{fig:prism} (middle, right). Since in an ideal tetrahedron the
  dihedral angles at opposite edges are equal~\cite{milnor82}
  \cite{milnor94:_volume_in_hyp}, most of the dihedral angles of the three
  tetrahedra are equal to some angle of the prism; see Figure~\ref{fig:prism}
  (right). The remaining two dihedral angles, $\lambda$ and $\mu$, are
  obtained by considering how the angles of the prism are sums of angles
  of the tetrahedra:
  \begin{equation*}
    \lambda=\beta'-\beta=\frac{\pi-\alpha-\beta-\gamma}{2},
    \qquad
    \mu=\pi-\gamma'.
  \end{equation*}
  Now the volumes of the three tetrahedra are obtained from
  Equation~\eqref{eq:Milnor}. Take the sum and note that $\ML(\pi-x)=-\ML(x)$
  to obtain Equation~\eqref{eq:vol_prism}.
\end{proof}

\subsection{The truncated volume of a tetrahedron with one hyperideal and
  three ideal vertices}

\begin{figure}[tbp]
  \centering \input{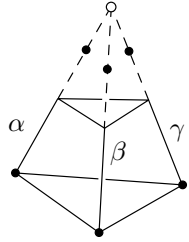}
  \caption{The tetrahedron $P_3$ with one hyperideal vertex (marked~$\circ$)
    and three ideal vertices. The points where the edges intersect the sphere
    at infinity are marked~$\bullet$. The tetrahedron is truncated at the
    polar plane of the hyperideal vertex.}
  \label{fig:tet_one_hyp}
\end{figure}

Let $P_3$ be a tetrahedron with one hyperideal and three ideal vertices as
shown in Figure~\ref{fig:tet_one_hyp}. 

\begin{lemma}
  The truncated volume (see Definition on p.~\pageref{pag:def_trunc_vol}) of
  $P_3$ is
  \begin{multline}
    \label{eq:vol_P3}
    \Vol(P_3) =
    \frac{1}{2}\Big(\ML(\alpha) + \ML(\beta) + \ML(\gamma)\\
    + \ML\big(\frac{\pi+\alpha-\beta-\gamma}{2}\big)
    + \ML\big(\frac{\pi-\alpha+\beta-\gamma}{2}\big)\\
    + \ML\big(\frac{\pi-\alpha-\beta+\gamma}{2}\big)
    + \ML\big(\frac{\pi-\alpha-\beta-\gamma}{2}\big)\Big).
  \end{multline}  
\end{lemma}

\begin{proof}
  If you reflect the truncated tetrahedron at the truncation plane, you get
  an ideal prism. The volume of the truncated tetrahedron is therefore half
  the volume of the ideal prism, Equation~\eqref{eq:vol_prism}.
\end{proof}

\begin{remark}
  The same volume formula holds when the apex is finite instead of
  hyperideal~\cite{vinberg91} \cite{vinberg93}. In that case
  $\alpha + \beta + \gamma > \pi$.
\end{remark}

\subsection{The truncated volume of a special pyramid}
\label{sec:spec_pyramid}

\begin{figure}[tbp]
  \centering \input{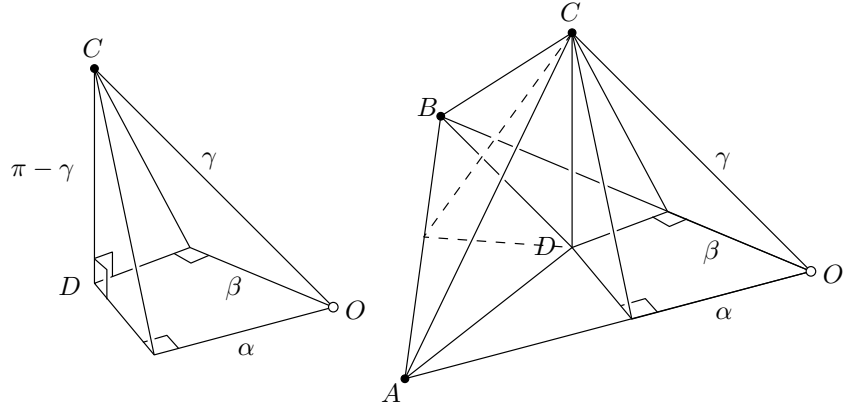}
  \caption{The special pyramid $P_4$ (left). Decomposition of a tetrahedron
    with one hyperideal and three ideal vertices into one special pyramid and
    four birectangular tetrahedra, two of which are mirror symmetric to each
    other (right).}
  \label{fig:special_pyramid}
\end{figure}

Let $P_4$ be a pyramid over a quadrilateral base, such that the apex $C$ is
at infinity, one lateral edge $CD$ is perpendicular to the base, the vertex
$O$ of the base that is opposite $D$ is hyperideal, and at the other two
vertices the angles of the base quadrilateral are $\frac{\pi}{2}$; see
Figure~\ref{fig:special_pyramid} (left). 
Let the interior dihedral angles at the edges emanating from $O$ be $\alpha$,
$\beta$, and $\gamma$ as shown. They satisfy $\alpha,\beta<\pi/2$ and
$\alpha+\beta+\gamma<\pi$.  Since the sum of dihedral angles at the
four-valent apex $C$ is $2\pi$ and two of the incident edges have dihedral
angle $\frac{\pi}{2}$, the dihedral angle at edge $CD$ is $\pi-\gamma$.

\begin{lemma}
  The truncated volume of $P_4$ is
  \begin{multline}
    \label{eq:vol_P4}
    V(P_4)=\frac{1}{2}\bigg(
    \ML(\gamma)
    + \ML\Big(\frac{\pi+\alpha-\beta-\gamma}{2}\Big)
    + \ML\Big(\frac{\pi-\alpha+\beta-\gamma}{2}\Big)\\
    - \ML\Big(\frac{\pi-\alpha-\beta+\gamma}{2}\Big)
    + \ML\Big(\frac{\pi-\alpha-\beta-\gamma}{2}\Big)
    \bigg).
  \end{multline}
\end{lemma}

\begin{proof}
  Extend the edges of the base emanating from $O$ until they intersect the
  infinite boundary at the ideal points $A$, $B$, see
  Figure~\ref{fig:special_pyramid} (right). The truncated volume of the
  tetrahedron $P_3$ with vertices $ABCO$ is given by
  Equation~\eqref{eq:vol_P3}. It can be partitioned into the special pyramid
  $P_4$, the tetrahedron $ABCD$, and two birectangular tetrahedra as shown.
  The volumes of the birectangular tetrahedra are given by
  Equation~\eqref{eq:vol_P1}; they are $\frac{1}{2}\ML(\alpha)$ and
  $\frac{1}{2}\ML(\beta)$. The tetrahedron $ABCD$ is symmetric with respect
  to reflection at the plane that contains edge $CD$ and intersects edge $AB$
  orthogonally. This symmetry plane splits the tetrahedron $ABCD$ into two
  symmetric birectangular tetrahedra. At edge $CD$, each of them has a
  dihedral angle of $\tfrac{1}{2}(\pi-\alpha-\beta+\gamma)$. (To see this
  consider the dihedral angles of all the pieces at edge $CD$; their sum is
  $2\pi$. Thus the angle in question is
  $\tfrac{1}{2}(2\pi-\alpha-\beta-(\pi-\gamma))$.)  The volume of the
  tetrahedron $ABCD$ is therefore
  $\ML(\tfrac{1}{2}(\pi-\alpha-\beta+\gamma))$. Subtract from the truncated
  volume of $P_4$ the volumes of tetrahedron $ABCD$ and the two birectangular
  tetrahedra to obtain Equation~\eqref{eq:vol_P4}.
\end{proof}

\begin{remark}
  The same volume formula holds when $O$ is finite instead of
  hyperideal~\cite{vinberg91} \cite{vinberg93}. In that case
  $\alpha+\beta+\gamma>\pi$.
\end{remark}

\begin{figure}[tbp]
  \centering 
  \input{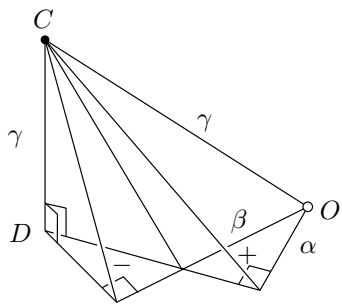}
  \caption{The base may also be a self-intersecting quadrilateral. Here,
    $\beta>\frac{\pi}{2}$.}
  \label{fig:special_pyramid_cross}
\end{figure}

Equation~\eqref{eq:vol_P4} holds also when the base is a self-intersecting
quadrilateral; see Figure~\ref{fig:special_pyramid_cross}.
In this case one of the angles $\alpha$ or $\beta$ is greater than
$\frac{\pi}{2}$. We have to regard the pyramid $P_4$ as a difference of two
tetrahedra (with signs as indicated in
Figure~\ref{fig:special_pyramid_cross}), and its volume as the difference of
the volumes of these tetrahedra. One can derive
Equation~\eqref{eq:vol_P4} by a similar construction.

\subsection{The truncated volume of a tetrahedron with one ideal and three
  hyperideal vertices}
\label{sec:vol_1_ideal_3_hyper}

\begin{figure}
  \centering
  \input{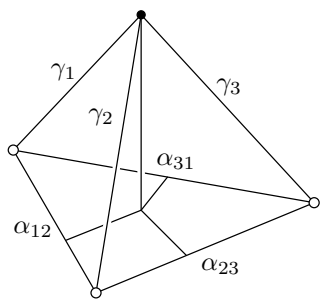}
  \caption{Tetrahedron with one ideal and three hyperideal vertices,
    partitioned into four special pyramids}
  \label{fig:four_pyramids}
\end{figure}

Finally, consider a tetrahedron with one ideal and three hyperideal vertices
as shown in Figure~\ref{fig:truncated_tet}. Drop the perpendicular from the
ideal vertex onto the plane of the opposite face and subdivide the
tetrahedron into four special pyramids (the bases of which may be
self-intersecting quadrilaterals) as shown Figure~\ref{fig:four_pyramids}.
The volume formula for the tetrahedron, Equation~\eqref{eq:V}, is obtained as
the sum of the volumes of the four special pyramids, given by
Equation~\eqref{eq:vol_P4}. (Note
$\ML(\frac{\pi}{2}-x)=-\ML(\frac{\pi}{2}+x)$.) This concludes the proof of
Lemma~\ref{lem:volume_V}.

\bibliographystyle{amsplain}
\bibliography{hyperideal}

\end{document}